\documentclass[3p,twocolumn]{elsarticle}

\usepackage{lineno,hyperref}
\modulolinenumbers[5]

\usepackage{graphicx}

\usepackage{mathptmx}      
\usepackage{amsmath} 
\DeclareMathOperator*{\argmin}{arg\,min}
\usepackage{amsfonts} 
\usepackage{amssymb}
\usepackage{subcaption}
\usepackage{xcolor}

\journal{my journal}

\begin{document}

\begin{frontmatter}

\title{Well control optimization using a two-step surrogate treatment}

\author{Daniel U. de Brito\fnref{myfootnote1,myfootnote2}}
\ead{dubrito@stanford.edu}
\author{Louis J. Durlofsky\fnref{myfootnote1}}
\ead{lou@stanford.edu}
\fntext[myfootnote1]{Department of Energy Resources Engineering, Stanford University, Stanford,CA 94305, USA}
\fntext[myfootnote2]{Now at Petr\'oleo Brasileiro S.A., Petrobras}




\begin{abstract}

Large numbers of flow simulations are typically required for the determination of optimal well settings. These simulations are often computationally demanding, which poses challenges for the optimizations. In this paper we present a new two-step surrogate treatment (ST) that reduces the computational expense associated with well control optimization. The method is applicable for oil production via waterflood, with well rates optimized at a single control period. The two-step ST entails two separate optimizations, which can both be performed very efficiently. In the first optimization, optimal well-rate ratios (i.e., the fraction of total injection or production associated with each well) are determined such that a measure of velocity variability over the field is minimized, leading to more uniform sweep. In the second step, overall injection and production rates are determined. The flow physics in the first step is highly simplified, while the actual physical system is simulated in the second step. Near-globally-optimal results can be determined in both cases, as the first optimization is posed as a QP problem, and the second step entails just a single optimization variable. Under full parallelization, the overall elapsed time for the ST corresponds to the runtime for 1--2 full-order simulations. Results are presented for multiple well configurations, for 2D and 3D channelized models, and comparisons with formal optimization procedures (mesh adaptive direct search or MADS, and an adjoint-gradient method) are conducted. Three different fluid mobility ratios ($M=1$, 3 and 5) are considered. Optimization results demonstrate that the two-step ST provides results in reasonable agreement with those from MADS and adjoint-gradient methods, with speedups of $5 \times$ or more. We also show that the ST is applicable in the inner-loop in field development optimization, where it will be especially useful since many different well configurations must be evaluated.

\end{abstract}


\end{frontmatter}



\section{Introduction}

\label{intro}

Well control optimization entails defining the set of controls (well rates or bottom-hole pressures) that maximize a predefined cost function, such as net present value or total oil production. Optimization computations generally rely on a reservoir simulator for the evaluation of proposed well control strategies. Since the flow models are typically time dependent and nonlinear, the solution of the resulting nonconvex optimization problem may be computationally demanding. This is even more of an issue when well locations and controls are optimized simultaneously. In such cases, the objective function evaluation for each set of well locations requires the determination of optimal well controls.  

In this study, we introduce a new two-step surrogate treatment (ST) for well control optimization for oil fields producing under waterflood. In the first step of the ST, we determine the well-rate ratios (i.e., the fraction of total injection or production associated with each well), for an appropriately defined unit-mobility-ratio problem. In the second step of the ST, we fix the well-rate ratios found in the first step, and then solve a full-physics optimization problem in one variable to determine the overall injection/production rate that maximizes net present value (NPV). Although the problems in the ST are idealized, they are posed such that we can efficiently find near-global optima in both steps (though this does not mean we determine the global optimum for the actual problem). 

Well control optimization problems have been addressed using both adjoint-gradient and derivative-free methods. Adjoint-gradient-based methods, which include those described in \cite{pallav,zhang,jansen}, are very efficient though they require access to reservoir simulation source code. Derivative-free optimization (DFO) methods represent a useful alternative when gradients are not available. Both stochastic and pattern search DFO methods have been applied in well control optimization problems. Stochastic methods used in this setting include genetic algorithms (GAs) and particle swarm optimization (PSO), which are both population-based evolutionary methods. These procedures were shown to be applicable for both well control~\cite{almeida,harding,foroud,Lushpeev} and well location optimization~\cite{emerick,Onwunalu}. Pattern-search or stencil-based optimization techniques, such as mesh adaptive direct search (MADS)~\cite{mads}, have also been used for well control optimization~\cite{Echeverria40,Echeverria42,masoud, ushmaev}. In this work both adjoint-gradient-based and MADS procedures will be compared to our new two-stage surrogate treatment.

The use of proxy or surrogate models is an effective way to reduce the cost of optimization in cases when the underlying full-physics model is expensive to simulate. Surrogate treatments have been applied for both well control and well placement optimization problems. For example, in an early study involving field development optimization, Rosenwald and Green~\cite{rosenwald_and_green} used single-phase flow models, for which linear programming techniques are applicable. Subsequent studies that used surrogate models in the form of statistical proxies, for well placement optimization, include those of Guyaguler et al.~\cite{Guyaguler}, Yeten et al.~\cite{yeten} and Farmer et al.~\cite{Farmer}. Reduced-physics surrogate models were used by Wilson and Durlofsky~\cite{wilson} for optimizing the development of shale gas reservoirs~\cite{wilson}, while Aliyev and Durlofsky~\cite{Aliyev2017} applied a sequence of upscaled models as surrogates for joint (well location and control) optimization in waterflooding problems. Coarse models were also used as surrogates in well placement optimization by Bukhamsin~\cite{Bukhamsin}.

Surrogate treatments for well control optimization, based on reduced-order modeling with proper orthogonal decomposition, have been developed by a number of researchers; see, e.g.,~\cite{doren,MarcoOpt,jincongl}, along with the review by Jansen and Durlofsky~\cite{Jansen2017}. Reduced-physics models based on streamline methods were shown to provide useful surrogates in water injection optimization~\cite{MarcoThieleATCE}. Well control optimization applications using streamlines have also been described by Park and Datta-Gupta~\cite{Park}, who used 3D streamlines to generate flood efficiency maps, and by Wen et al.~\cite{tailai}, who applied both streamline and time-of-flight (TOF) computations for the optimization of mature fields under waterflood. Flow diagnostics tools~\cite{Shavali,Lie}, typically constructed based on reduced-physics simulations, represent another effective set of surrogate treatments that can be used for optimization. For example, M$\o$yner et al.~\cite{Moyner} performed well control optimization for systems involving two- and three-phase flow. Lie et al.~\cite{Lie} used quantities derived from TOF and multiscale methods as fast proxies for well control optimization under waterflood. We note finally that He et al.~\cite{jincong_proxy} and Chen et al.~\cite{CHEN_proxy_hm} have demonstrated that surrogate models can also be applied to accelerate history matching and uncertainty quantification computations.

Our two-step ST somewhat resembles the procedure developed by Rodr{\'i}guez Torrado et al.~\cite{torrado} and Embid Droz et al.~\cite{patent} for setting  well controls. Their methodology also relies on the optimization of a single overall field rate variable (e.g., field oil production rate), and was used in combination with well scheduling and location optimization procedures. Well controls with this approach are determined through estimation of the field recovery factor, while the field rate variable is obtained using derivative-free methods. A key difference between our two-step ST and the procedure for optimizing the well controls developed by Rodr{\'i}guez Torrado et al.~\cite{torrado} is in the (heuristic) scheme used to determine the well-rate ratios. In our approach, a simplified flow problem is considered for the determination of the well-rate ratios, while their procedure uses static reservoir properties. Detailed comparisons would, however, need to be performed to clearly establish the relative advantages of the two methods.

In this paper, we develop and apply a new two-step surrogate treatment (ST) to enable efficient well control optimization for subsurface flow problems. Although field development optimization is not the focus of this paper, a key applications for our ST is the joint optimization of well locations and controls, where the ST can be used as a fast optimization method for the (inner-loop) well control optimization. The computational savings achieved by ST derive from solving two approximate subproblems rather than the actual full-physics problem. In the first subproblem we simplify the physics, and in the second subproblem we consider the actual physics but greatly simplify the optimization problem. Optimum solutions can be obtained for both subproblems at a cost of 1--2 full-order runs under full parallelization. The method provides well controls for only a single control period, but this solution can then be used as an initial guess for cases involving multiple control periods. Geological uncertainty is not considered, though this could be addressed by optimizing over multiple realizations.

This paper is organized as follows. In Sect.~\ref{sec:2} we present the basic oil-water flow equations and the general well control optimization problem. In Sect.~\ref{sec:3}, the two-step surrogate treatment is described in detail.  Heuristic procedures for handling BHP constraints are discussed, and the general workflow is explained. In Sect.~\ref{sec:4}, we present ST well control optimization results for 40 different (waterflood) well scenarios in a 2D model. Well control optimization for additional mobility ratios, and comparisons with optimization results using adjoint-gradient and MADS procedures, are also provided. Results for a 3D case are then presented. Finally, in Sect.~\ref{sec:5}, we summarize this work and provide suggestions for future research in this area.


\section{Problem statement}
\label{sec:2}

In this section, we present the equations governing reservoir flow, and we provide a formal statement of the general optimization problem.

\subsection{Oil-water flow equations}
\label{subsec:2.1}

For two-phase oil-water flow, the governing equations, which derive from Darcy's law for each phase and mass conservation statements, can be expressed in terms of the so-called pressure and saturation equations. In the absence of capillary pressure and gravity effects, these equations are as follows:
\begin{equation}\label{pressure_eqn}
	\nabla\cdot \left(\lambda_t\left(S_w\right) \; \textbf{k} \nabla p\right)=-\tilde{q}^{w}_{t}+\phi c_t \frac {\partial p} {\partial t} ,
\end{equation}

\begin{equation}\label{saturation_eqn}
	\phi \frac{\partial S_w}{\partial t} +\nabla \cdot\left(f\left(S_w\right) \; \textbf{u}_t \right)=\tilde{q}^{w}_{w}.
\end{equation}
\smallskip
\noindent Here subscripts $o$ and $w$ refer to oil and water phases, $\lambda_t$ is the total mobility, with $\lambda_t=\lambda_o+\lambda_w$, where $\lambda_o=k_{ro}/\mu_o$ and $\lambda_w=k_{rw}/\mu_w$, $k_{ro}$ and $k_{rw}$ are relative permeability to the oil and water phases respectively, and $\mu_o$ and $\mu_w$ are phase viscosities. Additional variables are water saturation $S_w$, absolute permeability tensor \textbf{k}, pressure $p$, total source/sink $\tilde{q}^{w}_{t}$ (the superscript $w$ here refers to wells), porosity $\phi$,  total compressibility
$c_t$ (which is a weighted combination of the fluid and rock compressibilities), and time $t$. The fractional flow of water $f(S_w)$ is defined as $f(S_w)=\lambda_w/\lambda_t$, the total velocity as $\textbf{u}_t=\textbf{u}_o+\textbf{u}_w$ ($\textbf{u}_o$ and $\textbf{u}_w$ are the Darcy velocities for oil and water), and $\tilde{q}^{w}_{w}$ is the water source/sink term. We additionally have the saturation constraint $S_o+S_w=1$, which completes the problem statement. 

In the discretized system, the source terms ($\tilde{q}^{w}_{w}$ and $\tilde{q}^{w}_{t}$) are represented using the Peaceman well model. These terms are nonzero only in grid blocks containing wells that are open to flow. The well model relates flow rate to bottom-hole pressure (BHP), and either well rates or BHPs can be specified as the control variables. In this study we consider vertical wells, completed over the entire reservoir thickness, with flow rates as the control parameter. Stanford's Automatic Differentiation-based General Purpose Research Simulator (AD-GPRS)~\cite{Zhou_ADGPRS} is used for all flow simulations.  Although gravitational terms are not shown in Eqs.~\ref{pressure_eqn} and \ref{saturation_eqn}, in order to simplify the presentation, these effects are included in the 3D simulations performed in Sect.~\ref{sec:4}. 

\subsection{General optimization problem}
\label{subsec:2.2}

The general field development optimization problem involves determination of the well types, locations and controls, with the goal of minimizing a cost function $J$. Following \cite{Isebor20143}, the optimization problem can be stated as follows:
\begin{gather}
\begin{array}{rrclcl}
\displaystyle \min_{\textbf{x} \in \mathbb{X},\textbf{u} \in \mathbb{U},\textbf{z} \in \mathbb{Z}} & {J (\textbf{x},\textbf{u},\textbf{z})}, \ \  \textrm{subject to} 
\begin{cases}
\textbf{g}(\textbf{p},\textbf{x},\textbf{u},\textbf{z})=\textbf{0},\\
\textbf{c}(\textbf{p},\textbf{x},\textbf{u},\textbf{z})\leq\textbf{0}
\end{cases}
\end{array}
\label{gen_field_dev_opt_eqn}
\end{gather}
\noindent The vectors \textbf{x} and \textbf{u} indicate integer (grid-block based) well location variables and continuous well control variables, respectively, while \textbf{z} are categorical variables, which indicate whether the well is an injector ($z_k=-1$), a producer ($z_k=1$), or not drilled at all ($z_k=0$).  The well location variables can also be treated as real-valued, and this may be preferable in cases where wells are not centered in grid blocks (e.g., with deviated wells). Here $\textbf{g}=\textbf{0}$ denotes the flow simulation equations, \textbf{p} represents the solution unknowns, which in our system are the pressure and saturation in every grid block, and $\textbf{c}$ defines any nonlinear constraints. The spaces $\mathbb{X}$ and $\mathbb{U}$ are defined to include bound constraints, which can be expressed as $\textbf{x}_l\leq \textbf{x}\leq \textbf{x}_u$ and $\textbf{u}_l\leq \textbf{u}\leq \textbf{u}_u$, where subscripts $l$ and $u$ denote lower and upper bounds.

In this work the objective is to maximize net present value (NPV); i.e., we set $J=-{\textrm{NPV}}$, with NPV given by: 
\begin{gather}
	\textnormal{NPV(\textbf{x},\textbf{u},\textbf{z})}=\sum\limits_{k = 1}^{n_p} \sum\limits_{s = 1}^{n_s} \frac{\Delta t_s \left(p_{o} ~q^{o}_{k,s}(\textbf{x},\textbf{u})-c_{pw}~q^{pw}_{k,s}(\textbf{x},\textbf{u}) \right)}{( 1+d )^{ \frac{t_s}{365}}}- \nonumber \\
    \sum\limits_{k = 1}^{n_i} \sum\limits_{s = 1}^{n_s} \frac{\Delta t_s~c_{iw}~q^{iw}_{k,s}(\textbf{x},\textbf{u})}{( 1+d )^{ \frac{t_s}{365}}} - \sum\limits_{k = 1}^{n_w}  \frac{|z_k|~c_w}{(1+d)^{ \frac{t_k}{365}}}.
    \label{gen_field_dev_npv_eqn}
\end{gather} 
\noindent Here $n_i$ is the number of injection wells, $n_p$ is the number of production wells, $n_w=n_i+n_p$ is the total number of wells, $n_s$ is the number of simulation time steps, $t_s$ and $\Delta t_s$ are the time and time step size at time step $s$, and $d$ is the annual discount rate. The price of oil, the cost for handling produced water, and the cost of injected water are, respectively, $p_o$, $c_{pw}$ and $c_{iw}$. The oil and water production rates and the water injection rate, for well $k$ at time step $s$, are denoted ${q}^o_{k,s}$, ${q}^{pw}_{k,s}$ and ${q}^{iw}_{k,s}$. The variable $t_k$ represents the time at which well $k$ is drilled, and the per-well drilling cost is denoted by $c_w$. 
 For a particular well control optimization problem, the last term in Eq.~\ref{gen_field_dev_npv_eqn} does not vary since the well configuration, and thus the values of $|z_k|$, are fixed. We nonetheless include the well costs in the computation in order to render the resulting NPVs more indicative of those for an actual project. \color{black}

In the optimizations performed in this study, the well types and locations are fixed (though in one set of examples we extract the well configurations from a field development optimization run). Then, with $J=-$NPV(\textbf{u)}, Eq.~\ref{gen_field_dev_opt_eqn} can be expressed as: 
\begin{equation}\label{rate_control_opt_eqn}
\begin{array}{rrclcl}
\displaystyle \max_{\textbf{u} \in \mathbb{U}} & {\textnormal{NPV}(\textbf{u})}, \ \  \textrm{subject to} 
\begin{cases}
\textbf{g}(\textbf{p},\textbf{u})=\textbf{0},\\
\textbf{c}(\textbf{p},\textbf{u})\leq\textbf{0}
\end{cases}
\end{array}
\end{equation}
\noindent The space $\mathbb{U}$ again includes the bound constraints for the continuous well control variables, and \textbf{c} specifies any nonlinear constraints.

We note that, given the nature of the optimization problem, fixing the well locations and types represents a considerable simplification. This specification eliminates the integer well placement variables (these variables are integers because well locations are represented in terms of discrete grid blocks), as well as the categorical well type variables. If these variables are included in the optimization, the overall field development optimization problem (Eq.~\ref{gen_field_dev_opt_eqn}) is a mixed integer nonlinear programming (MINLP) problem, which is much more difficult to solve than the real-variable problem defined in Eq.~\ref{rate_control_opt_eqn}.

The MINLP field development problem has been addressed using both sequential and joint optimization procedures. In the sequential approach, the well types and locations are optimized first, with well controls, or a well control strategy, specified. Once the optimized well placement configuration is determined, the well controls are then optimized for this configuration. However, since the optimal well locations depend on how the wells are operated, the sequential approach will in general provide a suboptimal solution. In the joint optimization approach, well types, locations, and controls are optimized together. This can be accomplished either in a single loop~\cite{Isebor20141}, in which all variables are considered simultaneously, or in a nested fashion~\cite{Bellout2012}. In the nested approach, well types and locations are defined in the outer loop, and well controls are optimized in the inner-loop. The ST developed in this study can be used directly in the inner-loop in nested (joint) field development optimization procedures. In this more challenging application, our two-step ST may provide significant computational speedup.

It is very useful for the surrogate model to be constructed such that the global optimum of the resulting optimization problem can be readily found. This will be the case, for example, if the surrogate treatment leads to a linear programming (LP) or a quadratic programming (QP) problem. For general LP or QP problems, a variety of global optimization algorithms are available. These include the Simplex algorithm and its extensions, and interior-point and conjugate-gradient methods~\cite{Dantzig1}.  As we will see, the two-step ST developed in this work can indeed be formulated as two supbroblems that are amenable to global optimization. We reiterate that this does not mean we find the global optimum for the original full-physics problem, but rather that we can find the global optima for the subproblems considered in the ST.

\section{Optimization methodology}
\label{sec:3}

In this section we present our two-step surrogate treatment for well control optimization of oil fields under waterflood. Both ST steps are described in detail, and the handling of BHP constraints is also discussed. This treatment considers only a single control period. Although our development here is for production via waterflood, ST is also expected to be applicable for other displacement processes, such as gas or water-alternating-gas injection, though it may be more approximate for these operations.

\subsection{Two-step ST} 
\label{subsec:3.1}

We denote the optimization variables used in this problem as $\textbf{f} \in \mathbb{R}^{n_w\times 1}$, where $n_w$ is the number of wells, and $Q \in \mathbb{R}$.  Here \textbf{f} represents the vector of well-rate ratios, or the fraction of the field rate that is allocated to each well, and \textit{Q} specifies the total field injection and production rates. At this point these total rates are the same since we have unit voidage replacement ratio, VRR = 1. \color{black}
With these definitions, the actual well rates can be obtained from the following relationship between \textbf{f} and \textit{Q}:
\begin{equation}
	\textbf{\^{q}}=\textbf{f}~Q,
    \label{ST_well_rates_eqn}
\end{equation}
\noindent where $\textbf{\^{q}}=\left[(\textbf{\^{q}\,\textsuperscript{\it{p}}})^{\textnormal{T}},(\textbf{\^{q}\,\textsuperscript{\it{i}}})^{\textnormal{T}}\right]^{\textnormal{T}} \in \mathbb{R}^{n_w\times 1}$, is the vector of well rates, $\textbf{f}=\left[(\textbf{f\,\textsuperscript{\it{p}}})^{\textnormal{T}},(\textbf{f\,\textsuperscript{\it{i}}})^{\textnormal{T}}\right]^{\textnormal{T}} \in \mathbb{R}^{n_w\times 1}$ and
$\sum\limits_{k = 1}^{n_p} {f}^{p}_{k}=\sum\limits_{k = 1}^{n_i} {f}^{i}_{k}=1$ (with all elements of \textbf{f} between 0 and 1), where superscripts $p$ and $i$ indicate production and injection wells. 

With this parameterization, the well control optimization problem can be expressed as:
\begin{equation}\label{gen_surrogate_eqn}
	\left[{\textbf{f}^{\ast}}^{\textnormal{T}},Q^{\ast}\right]^{\textnormal{T}}=\displaystyle \argmin_{\textbf{f} \in F, Q\in \mathbb{R}} {J (\textbf{f},Q)},
\end{equation}
\noindent where \textit{J} is the objective function of interest (negative NPV in our case) and \color{black} $F \in \mathbb{R}^{n_w \times 1} \subset [0,1]$. \color{black}


We address the optimization problem stated in Eq.~\ref{gen_surrogate_eqn} using a two-step optimization approach. This choice allows us to break the well control optimization problem into two subproblems involving different objective functions. In the first step, we optimize \textbf{f}:
\begin{equation}\label{gen_1st_step_surrogate_eqn}
	\textbf{f}^{\ast}=\displaystyle \argmin_{\textbf{f} \in F} {J_a (\textbf{f})},
\end{equation}
\noindent where $J_a$ is the first objective function. In the second step we optimize \textit{Q} using $\textbf{f}^\ast$ from the first step:
\begin{equation}\label{gen_2nd_step_surrogate_eqn}
	Q^{\ast}=\displaystyle \argmin_{Q\in \mathbb{R}} {J_b (\textbf{f}^{\ast},Q)},
\end{equation}
\noindent where $J_b$ is the second objective function.

The objective function used in the determination of optimal well-rate ratios corresponds to a measure of velocity variability. By minimizing this quantity, a more uniform sweep is achieved. This minimization is formulated as a QP problem, under the assumptions of unit mobility ratio ($M=1$) and VRR = 1. Given the optimal well-rate ratios $\textbf{f}^{\ast}$, $Q$ is determined in the second subproblem, by running a limited number of full-physics simulations, such that the NPV of the actual problem is maximized. We now describe the two subproblems in turn.


\subsubsection{First subproblem: optimization of \textnormal{\textbf{f}} }
\label{subsubsec:3.1.1}

We introduce significant approximations, in the first step of the ST, relative to the (target) optimization of the full-physics problem. Specifically, our goal in the first step is to minimize a measure of velocity variation within the model, which we expect to lead to improved sweep and thus higher oil production and NPV (this expectation will be confirmed below). In this step we assume unit mobility ratio, which results in a linear flow equation. We note that other simulation-based approaches such as flow diagnostics~\cite{Moyner,Lie} could also be used in this step, though it is not clear that such approaches would provide problems amenable to global optimization.

The $M=1$ pressure equation can be recovered from the two-phase flow system (Eqs.~\ref{pressure_eqn} and \ref{saturation_eqn}) by setting $\mu_w = \mu_o=\mu$, and the relative permeabilities equal to their corresponding phase saturations; i.e., $f_w(S_w)=k_{rw}=S_w$ and $f_o(S_o)=k_{ro}=S_o=1-S_w$. These assumptions simplify the governing equation, which in the case of slightly compressible systems can be expressed as:
\begin{equation}\label{sligthly_compress_systems_pressure_eqn}
	\nabla \cdot \left(\frac{\textbf{k}}{\mu} \nabla p\right)=-\tilde{q}^{w}_{t}+\phi c \frac{\partial p}{\partial t},
\end{equation}
\noindent where $c$ is compressibility. Darcy velocity in this case is given by:
\begin{equation}\label{total_vel_unit_mobility_eqn}
	\textbf{u}=-\frac{\textbf{k}}{\mu}\nabla p.
\end{equation}

In our implementation, we solve Eq.~\ref{sligthly_compress_systems_pressure_eqn} numerically using a standard finite volume implementation. Velocity components, constructed from this pressure solution, are used in the first-stage optimization. For a 3D Cartesian system, with blocks indexed as $i$, $j$ and $k$ in the $x$, $y$ and $z$-directions, and corresponding block interfaces designated $i\pm1/2$, $j\pm1/2$ and $k\pm1/2$, the velocity components in the $x$-direction, across block faces $i\pm1/2$, are given by: 
\begin{equation}\label{xvel_comp_eqn}
	({u}_x)_{i\pm\frac{1}{2},j,k}=\frac{(\tilde{q})_{i\pm\frac{1}{2},j,k}}{A}.
\end{equation}
\noindent Here $A$ is the area of the interface and $(\tilde{q})_{i\pm\frac{1}{2},j,k}$ are the block-to-block flow rates. These are given by, e.g., $(\tilde{q})_{i-\frac{1}{2},j,k}=(T_x)_{i-\frac{1}{2},j,k} \Delta p$, where transmissibility $(T_x)_{i-\frac{1}{2},j,k}=\frac{k_h A}{\mu \Delta x}$, where $k_h$ is the weighted harmonic mean of $k_x$ in blocks $i,j,k$ and $i-1,j,k$, $\Delta x$ is the center-to-center distance between these two blocks, and $\Delta p$ is the pressure difference between the two blocks. Analogous expressions provide the fluxes and velocity components across the other grid block faces.

The velocity components at grid block centers can be obtained by averaging the velocity components at the faces: 
\begin{equation}\label{xvel_comp_avg_eqn}
	({u}_x)_{i,j,k}=\frac{({u}_x)_{i-\frac{1}{2},j,k}+({u}_x)_{i+\frac{1}{2},j,k}}{2}.
\end{equation}
\noindent We denote the velocity components for grid block $i,j,k$ as $(u_d)_{i,j,k}$, with $d=x,y,z$. 

The linearity of Eq.~\ref{sligthly_compress_systems_pressure_eqn} implies that both the pressure and the velocity fields can be expressed as a superposition of the responses of individual wells~\cite{Ahlfeld,SudaryantoYortsos}. This also allows the use of response matrices for pressure and velocity computations~\cite{rosenwald_and_green,Ahlfeld}, which can accelerate considerably the calculation of such quantities when well rates are changed. Specifically, the velocity response matrices $\textbf{V}_d \in \mathbb{R}^{n_b\times n_w}$, where $n_b$ is the number of grid blocks, enable the resulting velocity fields $\textbf{\^{u}}_d \in \mathbb{R}^{n_b\times 1}$ to be written as linear combinations of the specified well rates $\textbf{\^{q}}$:
\begin{equation}\label{vel_field_eqn_x}
	\textbf{\^{u}}_d=\textbf{V}_d ~\textbf{\^{q}},
\end{equation}                                      
\noindent where $d=x,y,z$.

The methodology to build the response matrices consists of running a number of single-phase flow simulations to obtain a series of velocity coefficients. These coefficients are simply the velocity field components at every grid block of the model resulting from a single well producing at a unit flow rate. For well control optimization problems, where the well number, type, and locations are fixed, the number of single-phase flow simulations required is equal to the number of wells $n_w$. In each simulation, one of the wells operates individually, at a constant flow rate, until the reservoir pressure reaches pseudo-steady state (PSS); i.e., until $\frac{\partial p}{\partial t}$ becomes a constant. When this occurs, the velocity components are saved as columns in the velocity response matrices.   

It is also possible (and faster) to compute the PSS solution directly rather than through the time-integration of Eq.~\ref{sligthly_compress_systems_pressure_eqn}. The direct PSS solution can be computed by replacing the right-hand side of Eq.~\ref{sligthly_compress_systems_pressure_eqn} by $\phi K$, where $K$ is an arbitrary constant. We use time integration here since AD-GPRS is applied for these computations, and this simulator does not currently provide the PSS solution directly.  We note finally that there are other ways to construct the velocity response matrices. For example, since injection and production balance in our problem, we could use $n_w-1$ runs, each involving two of the wells (instead of the $n_w$ single-well runs used in this work). The construction of the response matrices is fast, however, and the additional PSS solution required by our approach has very little impact on total elapsed time. \color{black}

\begin{figure*}
  \begin{center}
  \includegraphics[width=0.80\textwidth]{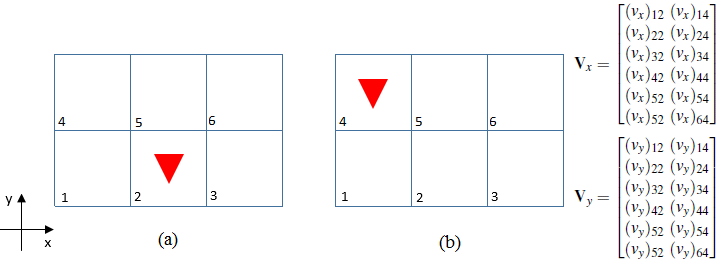}
  \end{center}
   \caption{Velocity response matrix for a $3 \times 2$ simulation model with two wells located at grid blocks 2 and 4.}
  \label{fig:vel_matrix_assemb2D}
\end{figure*}

Figure~\ref{fig:vel_matrix_assemb2D} depicts the assembly of the velocity response matrix  for a $3 \times 2$ model with wells located at grid blocks 2 and 4. The two model setups in Fig.~\ref{fig:vel_matrix_assemb2D}a,b indicate the well configurations for the two simulation runs required to build the response matrices. We denote $(v_d)_{ij}$, an element of $\textbf{V}_d$, to be the velocity component for grid block $i$ due to the operation of a single well located in grid block $j$. The process is analogous for 3D cases.   

\begin{figure}[!htb]
    \centering
    \includegraphics[width=0.5 \textwidth]{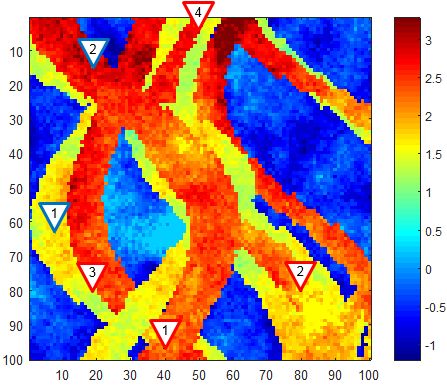}
    \caption{2D channelized model containing four production wells (shown in red) and two injection wells (blue). Permeability field ($\log_{10}k$, with $k$ in mD) generated by Isebor and Durlofsky~\cite{Isebor20143}. Table~\ref{tab:100x100_perm} summarizes other relevant simulation parameters.}
    \label{fig:100x100_4prod2inj_superposition}
\end{figure}

\begin{figure*}[t!]
    \centering
    \begin{subfigure}[b]{0.32\textwidth}
        \centering
        \includegraphics[width=\textwidth]{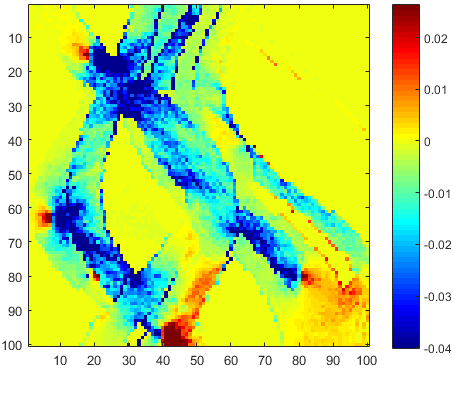}
        \caption{${\textbf{\^{u}}}^{\textnormal{super}}_x$}
    \end{subfigure}%
    ~ 
    \begin{subfigure}[b]{0.32\textwidth}
        \centering
         \includegraphics[width=\textwidth]{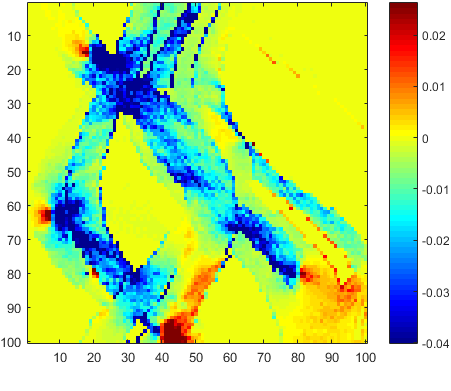}
        \caption{${\textbf{\^{u}}}^{\textnormal{simul}}_x$}
    \end{subfigure}
     ~ 
    \begin{subfigure}[b]{0.32\textwidth}
        \centering
        \includegraphics[width=\textwidth]{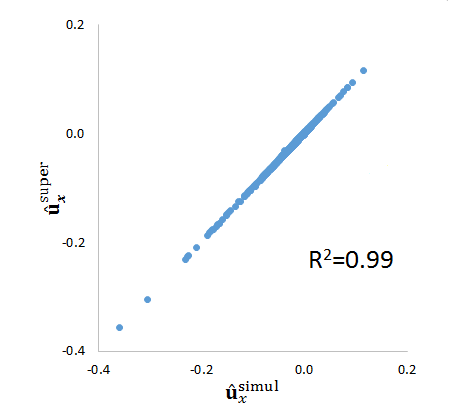}
        \caption{}
    \end{subfigure}
     ~ 
    \begin{subfigure}[b]{0.32\textwidth}
        \centering
        \includegraphics[width=\textwidth]{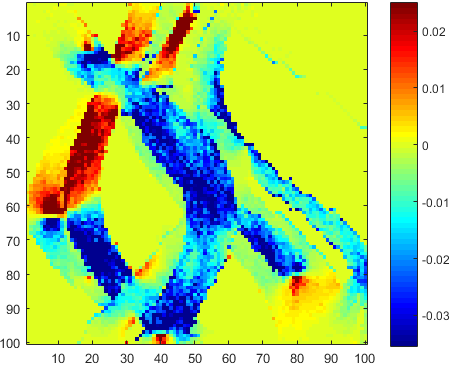}
        \caption{${\textbf{\^{u}}}^{\textnormal{super}}_y$}
    \end{subfigure}
     ~ 
    \begin{subfigure}[b]{0.32\textwidth}
        \centering
        \includegraphics[width=\textwidth]{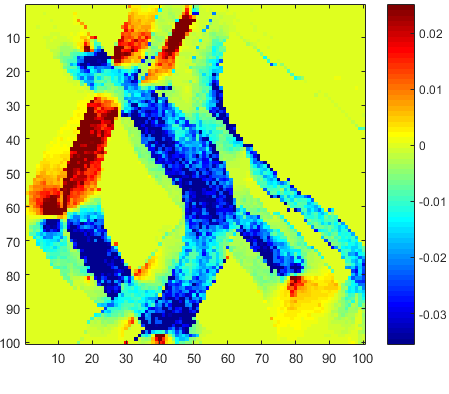}
        \caption{${\textbf{\^{u}}}^{\textnormal{simul}}_y$}
    \end{subfigure}
     ~ 
    \begin{subfigure}[b]{0.32\textwidth}
        \centering
        \includegraphics[width=\textwidth]{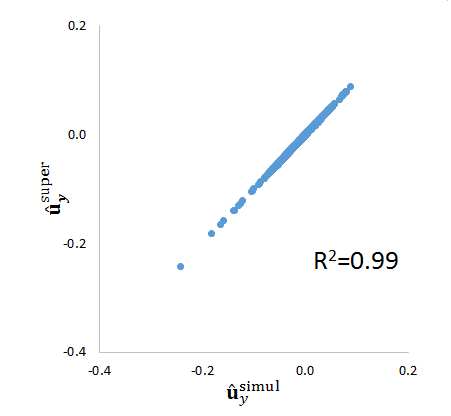}
        \caption{}
    \end{subfigure}   
    \caption{Velocity components (for model shown in Fig.~\ref{fig:100x100_4prod2inj_superposition}) obtained through superposition and direct simulation. (a,b,c) depict $x$-direction components and cross-plot, and (d,e,f) depict $y$-direction components and cross-plot.}
    \label{fig:100x100_super_vs_simul_vel}
\end{figure*}

Even though the $\textbf{V}_d$ matrices are constructed from pseudo-steady state solutions considering one well at a time, the multiwell response obtained through superposition will reproduce the steady-state solution as long as $\textnormal{VRR}=1$. This is now demonstrated by comparing the velocity fields generated through direct simulation (i.e., by solution of Eq.~\ref{sligthly_compress_systems_pressure_eqn} until steady state is reached) and by superposition for a $100 \times 100$ model containing four production and two injection wells. The channelized model and well locations are shown in Fig.~\ref{fig:100x100_4prod2inj_superposition}. The resulting velocity fields in the $x$ and $y$-directions, as well as cross plots comparing the two solutions, are presented in Fig.~\ref{fig:100x100_super_vs_simul_vel}. The velocity fields constructed via direct simulation and by superposition agree closely. The average absolute errors for velocities computed through superposition are about $7\times10^{-5}$ for both components. Although the concept and use of velocity response matrices is not new~\cite{Ahlfeld}, we are not aware of this approach being used before to construct steady-state solutions from the superposition of individual well responses obtained at pseudo-steady state. 

The fact that the velocity field can be expressed as a linear combination of the well rates, using Eq.~\ref{vel_field_eqn_x}, means that the minimization of the (weighted) squared velocities can be formulated as a QP problem. As indicated earlier, through use of such a minimization, we can find the set of well rates $\textbf{f}$ that result in a more uniform sweep and thus (expected) higher oil recovery and/or NPV. The general QP minimization can be stated as follows:
\begin{gather}
\displaystyle\min_{\textbf{f}\in F} \ \ {{\textbf{\^{u}}}^{\textnormal{T}}_x \textbf{W}_x {\textbf{\^{u}}_x} + {\textbf{\^{u}}}^{\textnormal{T}}_y \textbf{W}_y {\textbf{\^{u}}_y} + {\textbf{\^{u}}}^{\textnormal{T}}_z \textbf{W}_z {\textbf{\^{u}}_z}},\nonumber \\
\text{subject to} \ \ \sum\limits_{k=1}^{n_p} {f}^p_k = \sum\limits_{k=1}^{n_i} {f}^i_k = 1, \ \ {f}^p_k,{f}^i_k \geq 0, \, f \in \mathbb{R},
\label{well_rate_ratios_opt_eqn}
\end{gather}
\noindent where $\textbf{W}_d \in \mathbb{R}^{n_b\times n_b}$ with $d=x,y,z$ is a diagonal matrix whose elements are positive weighting factors. The field rate $Q$ is arbitrary at this point, and here we set $Q=1$. In the results below we specify $\textbf{W}_d=\textbf{I}$ for simplicity. However, by setting the elements of $\textbf{W}_d$ more generally, we can emphasize velocity uniformity in particular areas of interest in the reservoir. This could be beneficial, for example, when new wells are introduced into a reservoir that has already been under production for some time.

Solution of the optimization problem in Eq.~\ref{well_rate_ratios_opt_eqn} provides well-rate ratios that minimize the sum of the resulting squared velocities over the entire model domain.  As we will see in Sect.~\ref{subsec:4.1.1}, this quantity does indeed correlate with NPV for the example considered. Furthermore, the QP problem from Eq.~\ref{well_rate_ratios_opt_eqn} involves a quadratic objective function subject to linear constraints, and it can be solved efficiently using well-established optimization techniques that guarantee that a global optimum is found (as discussed in Sect.~\ref{subsec:2.2}). The use of the response matrix approach enables the very fast computation of ${\textbf{\^{u}}}_d$ (i.e., we do not need to solve any flow equations), thus accelerating the determination of $\textbf{f}^{\ast}$. Note that the minimization of the actual velocity variance, rather than the squared velocity, also corresponds to a QP problem. We implemented both of these minimizations, and found that the minimization of squared velocity provided better results for the cases evaluated.  This may be because the use of squared velocity penalizes high velocities more strongly. Because high extremes in the velocity field can lead to early breakthrough and inefficient displacement, and thus suboptimal NPVs, the avoidance of these effects may be very beneficial.
\color{black}

The computational cost for constructing the response matrix and solving the optimization problem in Eq.~\ref{well_rate_ratios_opt_eqn} is around 30--40\% of a single full-physics AD-GPRS simulation for the cases considered here. However, as noted earlier, the PSS problem could be solved directly, without integrating in time. Such a treatment, though not implemented in this work, will decrease the cost of this step to just a few percent of a full-physics AD-GPRS run. 

 We reiterate that the optimization problem posed in this first ST step is heuristic. The problem does not include many of the quantities that contribute directly to NPV, such as oil price, water costs and discount rate. In addition, this problem entails simplified flow physics. Thus it is useful to assess the correspondence between minimizing squared velocity (Eq.~\ref{well_rate_ratios_opt_eqn}) and maximizing NPV. This will be considered in Sect.~\ref{subsec:4.1.1}.


\subsubsection{Second subproblem: optimization of $Q$}
\label{subsubsec:3.1.2}

Once the optimal well-rate ratios $\textbf{f}^{\ast}$ have been determined, the actual well rates become a linear function of the field rate $Q$, as indicated by Eq.~\ref{ST_well_rates_eqn}. Thus, the optimization now involves only a single variable. At this stage we consider the actual objective function (in this case NPV as defined in Eq.~\ref{gen_field_dev_npv_eqn}) and we simulate the full-physics system.

\begin{figure}[!htb]
    \centering
    \includegraphics[width=0.47 \textwidth]{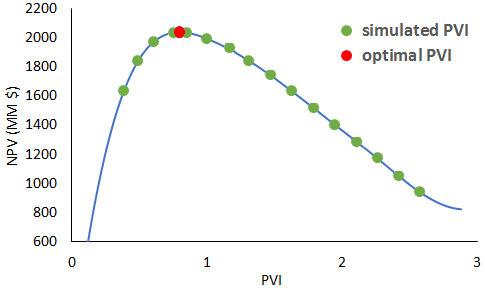}
    \caption{Relationship between NPV and PVI, which is used to determine the optimal field rate $Q^*$ in the second subproblem. The green points indicate results for different simulation runs, and the red point indicates the true optimum.}
    \label{fig:PVIvsNPV_curve_nonunitmobility}
\end{figure}

The second-stage optimization problem can be written as:
\begin{gather}
\displaystyle\max_{Q\in \mathbb{R}} \ \  \textnormal{NPV}(\textbf{f}^{\ast},Q), \ \ \text{subject to} \ \ Q > 0.
\label{field_rate_opt_eqn}
\end{gather}
At this stage we do not require $\textnormal{VRR}=1$ or the assumption of slight compressibility, though it is expected that the overall ST will perform better when these conditions are approximated. Since the well-rate ratios $\textbf{f}^{\ast}$ are specified, $\textnormal{VRR}=1$ will still be naturally enforced unless other constraints, such as well BHP limits (discussed below), become active. We treat the problem defined in Eq.~\ref{field_rate_opt_eqn} very simply. Specifically, we define a range for $Q$ such that, over the simulation time frame, between 0.5~pore volumes injected (PVI) and 2.5~PVI are introduced into the model from all of the wells. We compute even increments in PVI based on $n_{\textrm{proc}}$, where $n_{\textrm{proc}}$ is the number of available processors  (or the number of available simulator licenses if we have fewer licenses than processors), and then simulate each $Q$ value on a different processor. For example, for $n_{\textrm{proc}}=41$, we simulate cases corresponding to $0.5, 0.55, 0.60, \ldots, 2.45, 2.5$~PVI. This provides the green points shown in Fig.~\ref{fig:PVIvsNPV_curve_nonunitmobility} in an elapsed (wall-clock) time corresponding to the time required for one full-physics simulation. We then take the $Q$ that provides the best NPV to be $Q^*$. 

We could of course iterate this procedure further, by running additional cases in the most promising PVI regions, to provide the global optimum (red point in Fig.~\ref{fig:PVIvsNPV_curve_nonunitmobility}). Given the level of approximation inherent in the overall surrogate treatment, however, this additional iteration may not be necessary, assuming a sufficient number of processors (and simulator licenses) are available.

We note finally that in the special case where the system actually corresponds to $M=1$ (or is very close to $M=1$), it is possible to use highly efficient approaches such as the time-of-flight computations implemented in the Matlab Reservoir Simulation Toolbox, MRST ~\cite{MRSTbook} in place of full-physics simulations. In such cases it will be very inexpensive to perform an additional iteration in the second stage of the ST.

At this point it is useful to review the treatments and approximations introduced into our two-stage ST. In the first stage, we minimize a measure of velocity variability under the assumption of $M=1$ and $\textnormal{VRR}=1$. The objective function used in this minimization is not the actual objective function, though it does correspond to a quantity that is expected to impact displacement efficiency. Then, given the optimized well-rate ratios from this first step, in the second step we maximize the actual objective function (NPV) using the full-physics model. Global optima can be determined for both problems -- the first-stage optimization is a QP problem, and the second-stage optimization involves only a single variable. Since the second-stage problem can be solved in an elapsed time corresponding to a single full-physics run, and given that the first-stage problem can be solved even faster, the elapsed time for an optimization using our ST is the time required for 1--2 full-physics runs. 

 Another aspect of our two-step ST is that there is not any iteration between the first and second steps; i.e., each problem is solved only once. In many two-step optimization procedures, there is iteration between the two steps. In such cases, convergence to a local minimum is often assured. The lack of iteration renders our approach very efficient, but we cannot claim that it convergences to a local minimum of the actual problem. In Sect.~\ref{subsec:3.2}, we will describe a treatment for BHP constraint handling that does entail iteration between the two steps, but this procedure would only be used if a BHP constraint violation is detected. We note finally that the settings provided by the two-step ST can be used as an initial guess for a traditional optimizer (e.g., MADS), and in this optimization we do expect to converge to a local minimum of the actual problem.
\color{black}

\subsection{BHP constraint handling techniques}
\label{subsec:3.2}

The two-step ST described above did not include the treatment of any constraints. In the case of well-rate specifications, as are applied here, the most important constraint is expected to be a BHP constraint.  We now describe two heuristic approaches for handling this constraint. The first approach (which is the one actually implemented in this study) is applicable when the constraint violations involve a relatively small fraction of wells. The second approach, which is more complicated and involves iteration between the two ST steps, might be required if a large fraction of the wells are under BHP control. \color{black}

We are concerned mainly with cases for which $\mu_o \gtrsim \mu_w$, and thus $M \gtrsim 1$. In such cases, less viscous (more mobile) fluid is injected into the reservoir, so the differences in pressure between injectors and producers tend to decrease with time. Therefore, if BHP constraint violations do occur, this will often be at the start of the simulation (we check for BHP constraint violations just after the initial transient period, which lasts for a few time steps). These constraint violations will be detected in the second stage of the ST, where full-physics simulations are performed to determine $Q^*$. In such cases, we specify the BHP for the affected well(s), at the maximum or minimum, and then reapportion the remaining injection and production. Specifically, well-rate ratios in this case are assigned as:
\begin{equation}\label{updated_well_rate_ratios_eqn_BHP}
	(f')^p_{k}=\frac{{f}^p_k}{\sum\limits_{k=1}^{n_{pr}} {f}^p_k}, 
\end{equation}      
\noindent where $(f')^p_{k}$ is the updated well-rate ratio for production well $k$, and $n_{pr}$ is the number of production wells that operate under rate control. An analogous procedure is applied to the injection wells when a BHP constraint violation occurs at the start of the simulation.

\begin{figure}[!htb]
    \centering
    \includegraphics[width=0.5 \textwidth]{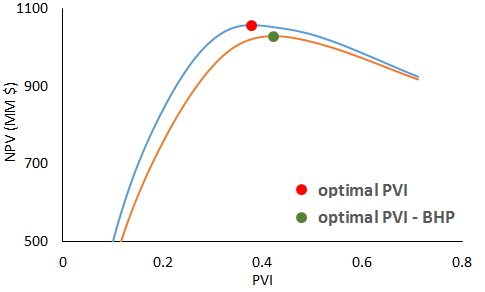}
    \caption{Relationship between NPV and PVI, used to determine the optimal field rate $Q^*$ in the second subproblem, for unconstrained (upper curve) and BHP-constrained (lower curve) cases. The points denote the optima.}
    \label{fig:PVIvsNPV_BHP_curve}
\end{figure}

The well-rate distributions now correspond to rescaled $\textbf{f}^\ast$ for those wells that remain under rate control. BHP constraints are satisfied for the other wells. We again vary rate $Q$ and compute NPV in the second ST stage. The resulting curve, which differs from that of the unconstrained case, is shown in Fig.~\ref{fig:PVIvsNPV_BHP_curve}. Note that the optimal rate $Q^\ast$ calculated in this step is not the actual total rate of the field. Rather, it is the total rate of the wells operating under rate control. The total rate of the field is the sum of $Q^\ast$ and the rate injected or produced by the wells operating under BHP control. On the other hand, this treatment does not impact the efficiency of the second-stage optimization -- it can still be achieved, in an elapsed time of 1--2 full-order simulations, as described above.

 When Eq.~\ref{updated_well_rate_ratios_eqn_BHP} is applied, the total rate of the field will not, in general, be strictly constant as a function of time. Field injection and production rates do, however, quickly `equilibrate,' meaning VRR=1 is again achieved (and sustained throughout the simulation). This typically occurs because another well also reaches a BHP constraint, and the specification of BHP for two wells, and rates for the other wells, results in balanced injection and production in a nearly incompressible system. Nonetheless, the $\textbf{f}^\ast$ values determined in the first ST step for the wells under BHP control are no longer honored. This acts to reduce the benefit afforded by minimizing the squared velocity in the first ST step. However, in cases where the fraction of BHP-controlled wells is relatively small, the overall impact on the two-step ST result is not expected to be large. We reiterate that the optimal rate determined in the second step fully accounts for wells operating under BHP control.  

BHP violations may also occur later in the simulation run  (i.e., after the early transient period, which lasts for a few time steps), though this is less of an issue for $M \gtrsim 1$, as noted earlier. If this happens, the well is switched to BHP control, and it remains on BHP control for the rest of the run. No additional corrective measures are taken, and the other wells are kept on rate control, unless they reach their BHP limit.

 Wells may also be specified to inject or produce under BHP control for operational reasons. In addition, under cases with high discount rate, early production will be preferred so rates may be higher. This could potentially lead to a larger fraction of wells reaching their BHP limits. Thus a more formal approach for handling BHP constraints may be required. We now describe such a treatment. Note that we have not implemented this approach since BHP constraints are not often reached in our runs, and the heuristic treatment described above suffices for present purposes.

 For cases where a large fraction of wells reach their BHP limits, or are specified to be under BHP control for operational reasons, we believe the following iterative procedure will be effective. In the first iteration, the two-step ST is performed exactly as described earlier. In the second ST step, however, after the early transient, we compute the maximum values of $f^i_k$ or $f^p_k$ corresponding to wells that reach their BHP limits. Then, the first ST step is repeated, with the BHP-constrained $f^i_k$ or $f^p_k$ specified as upper bounds, for the corresponding wells, in the QP optimization. This will ensure that the rate allocation is indeed optimized (under the reduced-physics assumptions in the first ST step) for the BHP-constrained case. The second ST step is then repeated to determine $Q^\ast$ given the new $\textbf{f}^\ast$ values. This entire procedure can be repeated (iterated) as necessary to assure no wells violate BHP constraints for the final $\textbf{f}^\ast$ and $Q^\ast$. We note that this iterative treatment directly links the two ST steps, since the results from second step impact the first step.
\color{black}

 A closely related issue is the treatment of wells that, for practical reasons, do not operate under rate control for extended periods, as discussed by \color{black} Ushmaev et al.~\cite{ushmaev} \color{black}. The strategy for handling such wells would be (essentially) that described above for wells that reach their BHP limit. In the case of wells specified to operate under BHP control, we would again iterate through the two-step ST, with upper-bound constraints for $f^i_k$ or $f^p_k$ determined from the second ST step. In this way the resulting $\textbf{f}^\ast$ and $Q^\ast$ would account for the fact that one or more wells are specified to operate under BHP control. 
\color{black}
 
 Practical problems may also involve nonlinear output constraints. Examples of these include the specification of maximum well water-oil ratio, or minimum well oil production rate. In general, such constraints can be handled both within the optimizer (using, e.g., a filter method as described in \cite{Isebor20141}), or directly in the simulator using well-control rules. We do not consider these types of constraints in this work. However, since the second ST step involves full-physics simulations, it would be straightforward to include these nonlinear constraints in the simulator. It might even be possible to account for them (approximately) in the iterative procedure described above.
\color{black}

\subsection{General workflow}
\label{subsec:3.3}

The two-step ST formulation entails the assembly of the velocity response matrices followed by the two optimization procedures. The Stanford Smart Fields Unified Optimization Framework (UOF)~\cite{yong,Sambandam} is used to test and compare the different optimization methods and treatments.

As explained earlier, for the construction of the velocity response matrices, a single-phase flow simulation is performed for each well, with the well specified to operate at constant flow rate. These simulations are run until the reservoir reaches pseudo-steady state. Once this is achieved, the pressure field is recorded, Darcy velocities are computed, and the velocity response matrices are assembled.

Algorithm~1 below describes the determination of ${\textbf f}^*$ and $Q^*$ using the two-stage ST. We first solve the squared-velocity-minimization problem (Eq.~\ref{well_rate_ratios_opt_eqn}) under the assumption of $M=1$ and $\textnormal{VRR}=1$. Then, NPV for the full-physics problem is maximized (Eq.~\ref{field_rate_opt_eqn}). The QP problem in the first step is solved in Matlab \cite{matlab} using CVX, a package for solving convex optimization problems~\cite{cvx}. The second optimization step requires running $n_{\textrm{proc}}$ full-physics simulations, where $n_{\textrm{proc}}$ is the number of computing nodes or simulator licenses available. If the full-physics model corresponds to $M=1$, MRST can be used to accelerate the second-stage computations.  

\begin{table*}[!htb]
    \centering
     \begin{tabular}{ll}
		\hline
        \multicolumn{2}{c}{ \textbf{Algorithm~1.} Well-rate ratios and field rate optimization} \\
        \hline
		1 &  	Solve first-step optimization problem to provide $\textbf{f}^{\ast}$ (Matlab/CVX) (Eq.~\ref{well_rate_ratios_opt_eqn}) \\
        2 &  	Solve second-step optimization problem to provide ${Q}^{\ast}$ (Eq.~\ref{field_rate_opt_eqn}): \\
        3 &  	\textbf{for} number of computing nodes (or licences) available \textbf{do} \\
        4 & 	\hspace{0.5cm} Run full-physics simulation with different field rates (AD-GPRS) \\
        5 &     \hspace{0.5cm} Treat BHP constraints as described in Sect.~\ref{subsec:3.2} \\
        6 &  	\textbf{end for} \\
        7 & 	Choose field rate associated with largest NPV\\    
        \hline
	\end{tabular}
    \label{tab:algorithm 2}
\end{table*}

\section{Example cases}
\label{sec:4}

\begin{figure}[!htb]
    \centering
    \includegraphics[width=0.5 \textwidth]{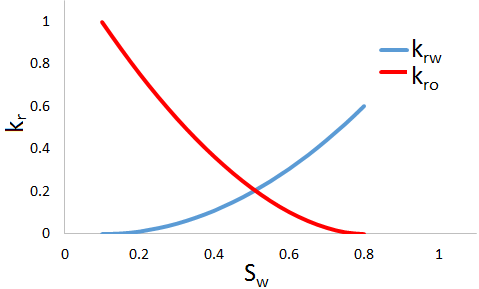}
    \caption{Relative permeability curves used for $M=3$ and $M=5$ simulations. Straight-line relative permeabilities are used for $M=1$ cases.}
    \label{fig:krel_curve}
\end{figure}

In this section, we apply the two-step ST described in Sect.~\ref{sec:3} to optimize well controls for 2D and 3D reservoir models. Performance is assessed for cases involving both single and multiple control periods. In the latter case, the single control period solution is used as the initial guess. The examples consider $M=1$, 3 and 5, where $M=\mu_o/\mu_w$. For $M$ = 3 and 5, we use the relative permeability curves shown in Fig.~\ref{fig:krel_curve}. For $M=1$, we specify $k_{rw}=S_w$, $k_{ro}=S_{o}$ and $\mu_w=\mu_o=1$. Thus the pressure and velocity fields for $M=1$ cases correspond exactly to those for single-phase flow. BHP constraints are imposed in these examples, but they are not violated at the start of the simulation in the cases presented. As noted earlier, Stanford's Automatic Differentiation-based General Purpose Research Simulator (AD-GPRS)~\cite{Zhou_ADGPRS} is used for all simulations.  Because the two-step ST is noninvasive with respect to the flow simulator, it can be used with any simulator. \color{black}

 Conventional optimization, using adjoint-gradients with SNOPT~\cite{snopt} and pattern search with MADS~\cite{mads}, is also performed in this study to enable comparison with ST results. These procedures are only guaranteed to find local optima, with the result depending on the initial guess, though in many cases the different local optima that are found correspond to very similar objective function values. This is not always the case, however. For example, in a recent study involving compositional models with nonlinear constraints, Kim~\cite{yong} found local optima that corresponded to very different objective function values. In this study we consider oil-water problems subject to only bound constraints, though we do observe variation in the objective function values, in some cases, using conventional methods. For this reason, five runs with different initial guesses are performed for the adjoint-gradient and MADS optimizations whenever these methods are applied.
\color{black}

In this study, adjoint-gradient optimizations are terminated when the gradient satisfies a threshold of 0.001, and the optimization converges to a local optimum, or after a maximum of 15 major iterations have been performed. MADS is terminated if the stencil size decreases beyond 0.01 (where this value corresponds to the scaled controls), or after a maximum of 15 MADS iterations (for 2D models) or 20 MADS iterations (for 3D models) have been performed. Each MADS iteration requires $2n_{\textrm{opt}}$ function evaluations (full-physics simulations in our case), where $n_{\textrm{opt}}$ is the number of optimization variables~\cite{mads}. Alternative stopping criteria for MADS, such as those introduced by Boch~\cite{Boch}, could also be considered.

\subsection{2D reservoir model}
\label{subsec:4.1}

Our first example involves the 2D reservoir model shown in Fig.~\ref{fig:100x100_4prod2inj_superposition} (this model was generated by \cite{Isebor20143}). This model, of dimensions $100\times100$, represents a complex channel-levee shale fluvial system. Relevant simulation and optimization parameters are provided in Table~\ref{tab:100x100_perm}. In our first set of examples, only a single control period (with well rates as the control variables) is used, so we have $n_{\textrm{opt}}=n_w-1$. 

\begin{table}[!htb]
    \centering
    \caption{Simulation and optimization parameters for $100 \times 100$ channelized model}
    \begin{tabular}{ll}
		\hline
		Grid dimensions                     & 100  $\times$ 100  $\times$ 1 \\
        Grid cell dimensions                & 50 ft $\times$ 50 ft $\times$ 100 ft \\
		Initial pressure $p_i$, at datum    & 6000 psi at 6000 ft\\
		$\mu_o$ 					        & 1.0, 3.0 or 5.0 cp \\
		$\mu_w$                 			& 1.0 cp \\
        $B_o$ and $B_w$ at $p_i$            & 1.075 RB/STB \\
		Oil compressibility at $p_i$		& 10$^{-5}$ psi$^{-1}$ \\
		Water compressibility at $p_i$		& 10$^{-5}$ psi$^{-1}$ \\
        Rock compressibility at $p_i$		& 10$^{-9}$ psi$^{-1}$  \\
        Porosity						    & 0.25\\
        Simulation period                   & 7300 days \\
        $p_o$, $c_{pw}$ and $c_{iw}$        & \$60, \$5 and \$5/STB \\
        $c_w$								& \$20 MM \\
        Discount rate                       & 0$\%$ \\
        \hline
	\end{tabular}
    \label{tab:100x100_perm}
\end{table}

\subsubsection{Results for different well scenarios}
\label{subsec:4.1.1}

Although field development optimization is not the main focus of this work, we first illustrate the potential use of the two-step ST within this context. This enables us to assess ST performance for a wide range of well configurations. These configurations are those actually generated during the course of a field development optimization run. Specifically, we ran the particle swarm optimization (PSO) algorithm described in \cite{Isebor20143}, with a PSO swarm size of 40. The goal was to maximize NPV by optimizing the number, type, locations and controls of a set of wells, with $M=3$. We extracted the configurations (well locations and types) corresponding to the 40 particles at a relatively early stage of the optimization -- specifically at iteration 25 (an early iteration was considered to maintain diversity within the swarm). We then applied ST for each of the 40 configurations. ST results are compared to optimization using MADS for the same sets of wells.

\begin{figure}[!htb]
    \centering
    \includegraphics[width=0.5 \textwidth]{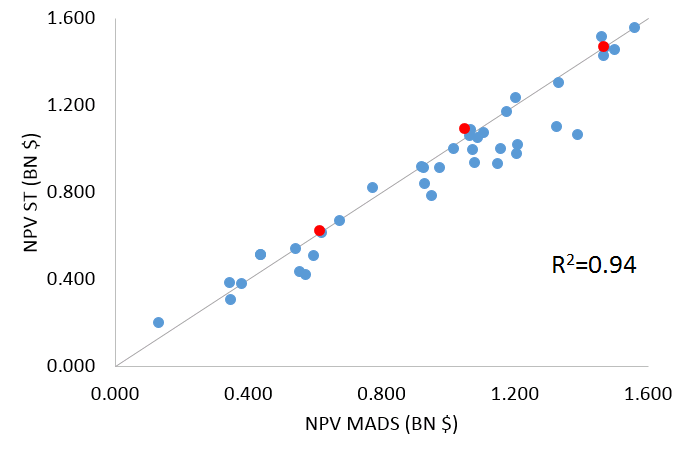}
    \caption{Comparison of NPVs obtained by ST and MADS for the 40 well scenarios extracted from a PSO field development optimization run ($M=3$). The red points indicate the three well scenarios used for further comparisons.}
    \label{fig:npv_st_vs_mads}
\end{figure}

Figure~\ref{fig:npv_st_vs_mads} displays a cross-plot of the NPVs obtained by ST and MADS. Each point represents one of the 40 PSO well scenarios. Relatively close correspondence between the ST and MADS NPVs is observed (correlation coefficient $R^2=0.94$), with most points falling near the 45-degree line.  These results demonstrate that the two-step ST is able to generate reasonable optimization results, over a range of NPVs, for the particular specifications in Table~\ref{tab:100x100_perm}.\color{black}~In other words, the two-step ST is able to perform satisfactorily for either poor, average or favorable sets of well types and locations. This is an important capability if the method is to be used for field development optimization, in which case a large fraction of the proposed well locations are far from optimal, and it is important that they be assessed accurately.

\begin{figure*}[t!]
    \centering
    \begin{subfigure}[b]{0.32\textwidth}
        \centering
        \includegraphics[width=\textwidth]{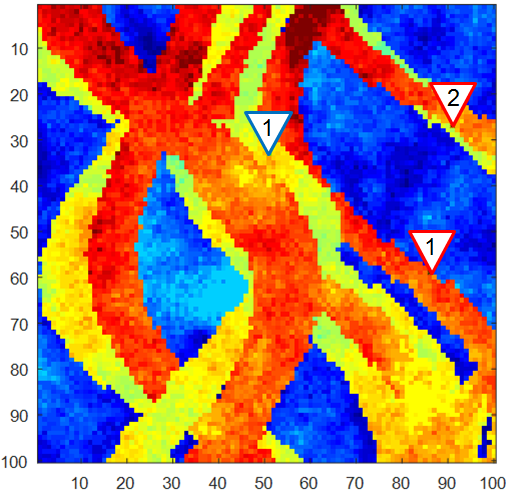}
        \caption{Well scenario 1}
    \end{subfigure}%
    ~ 
    \begin{subfigure}[b]{0.32\textwidth}
        \centering
        \includegraphics[width=\textwidth]{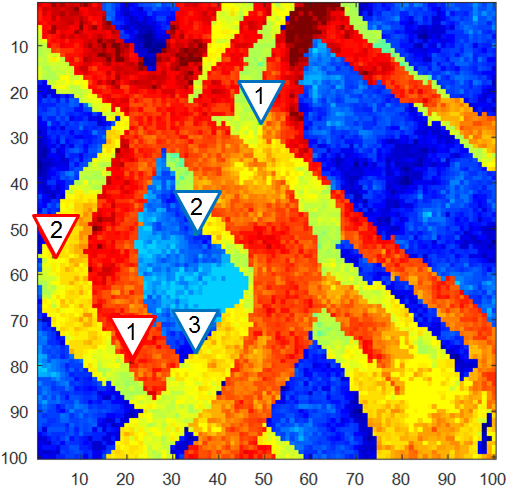}
        \caption{Well scenario 2}
    \end{subfigure}
     ~ 
    \begin{subfigure}[b]{0.32\textwidth}
        \centering
        \includegraphics[width=\textwidth]{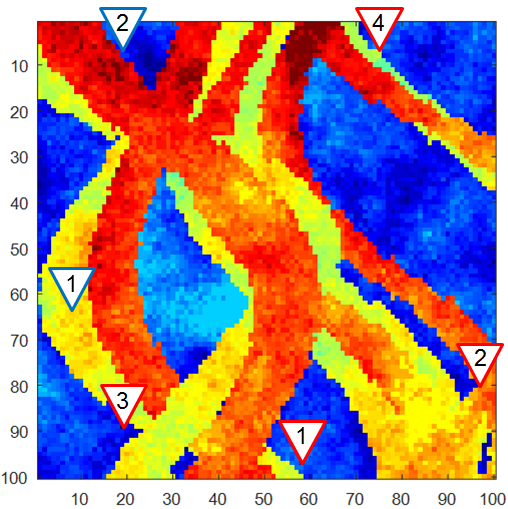}
        \caption{Well scenario 3}
    \end{subfigure}
     ~ 
    \begin{subfigure}[b]{0.32\textwidth}
        \centering
        \includegraphics[width=\textwidth]{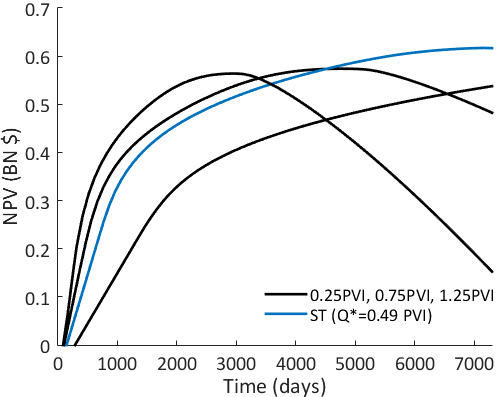}
        \caption{NPV profile for well scenario 1}
    \end{subfigure}
     ~ 
    \begin{subfigure}[b]{0.32\textwidth}
        \centering
        \includegraphics[width=\textwidth]{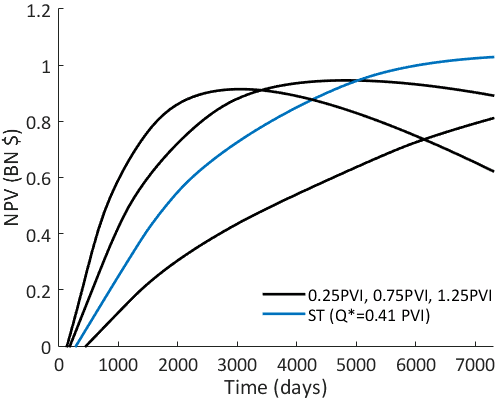}
        \caption{NPV profile for well scenario 2}
    \end{subfigure}
     ~ 
    \begin{subfigure}[b]{0.32\textwidth}
        \centering
        \includegraphics[width=\textwidth]{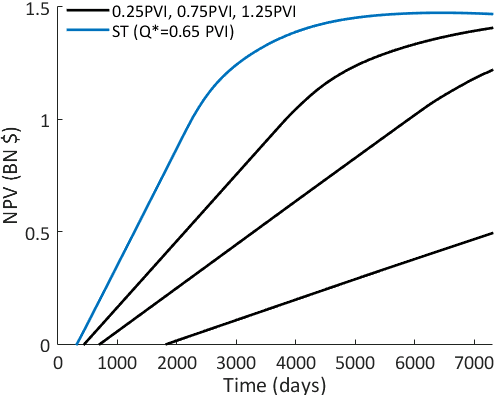}
        \caption{NPV profile for well scenario 3}
    \end{subfigure}   
    \caption{Well configurations and permeability field, and NPV evolution in time, for the three cases indicated by the red points in Fig.~\ref{fig:npv_st_vs_mads}. In (a, b, c), production wells are shown in red and injection wells in blue. Curves in (d, e, f) correspond to ST results (blue curves) and base-case results for 0.25~PVI, 0.75~PVI and 1.25~PVI. $M=3$ in all cases.}
    \label{fig:100x100_3example_cases}
\end{figure*}

We now perform further assessments for some of the well configurations considered in Fig.~\ref{fig:npv_st_vs_mads}. These scenarios are represented by the red points in Fig.~\ref{fig:npv_st_vs_mads}, and the corresponding well configurations are shown in Fig.~\ref{fig:100x100_3example_cases}a,b,c. It is evident from Fig.~\ref{fig:npv_st_vs_mads} that these cases cover a broad NPV range; i.e., we do not consider only the most promising well configurations.

 In the first assessment, well scenario~3 (Fig.~\ref{fig:100x100_3example_cases}c) is considered. Our intent here is to evaluate the correlation between NPV and the sum of the squared velocities over the entire domain (this is the quantity minimized in the first ST step). In this evaluation, we specify $M=1$ and VRR~=~1. In addition to the optimal well-rate ratios computed in the first ST step (${\textbf f}^*$), we consider 39 additional \textbf{f} vectors. One of these vectors corresponds to equally distributed rates, and the other 38 are generated randomly. In all cases the $f_k^p$ and $f_k^i$ sum to unity. For each of the 40 \textbf{f} vectors, the optimal field rate (${Q}^{*}$) is computed using the second ST step. This ${Q}^{*}$ determination is performed (separately) for discount rates $d$ of both 0 and 0.1 (10\%).

\begin{figure}[h]
    \centering
    \begin{subfigure}[b]{0.50\textwidth}
        \centering
        \includegraphics[width=\textwidth]{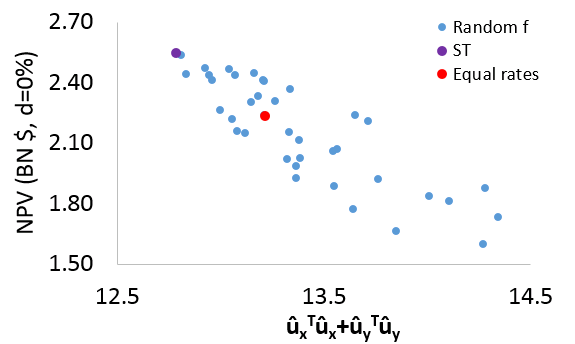}
        \caption{0\% annual discounting}
    \end{subfigure}%
   \bigbreak
    \begin{subfigure}[b]{0.50\textwidth}
        \centering
        \includegraphics[width=\textwidth]{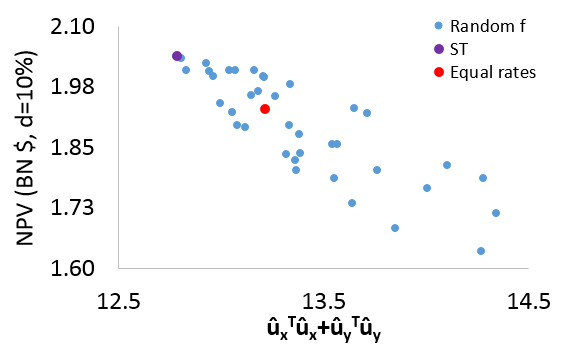}
        \caption{10\% annual discounting}
    \end{subfigure}  
    \caption{ Correspondence between NPV and sum of the squared velocities over the entire domain, for well scenario~3 (Fig.~\ref{fig:100x100_3example_cases}c) with $M=1$. NPVs computed for ST result (${\textbf f}^*$), equally distributed rates, and 38 random \textbf{f} vectors. Separate evaluations performed for (a) 0\% discount rate and (b) 10\% annual discount rate. \color{black}}
    \label{fig:npv_vs_uTu}
\end{figure}

Figure~\ref{fig:npv_vs_uTu} displays plots of NPV versus the corresponding sum of the squared velocities for each of the 40 \textbf{f} vectors, for $d=0$ and 0.1. The dark blue point in each plot corresponds to the two-step ST result, the red point to the equally distributed rates case, and the light blue points to the 38 cases with random \textbf{f} vectors. The scale for the $x$-axis is arbitrary, since the results of the first ST step do not depend on $Q$ (here we use $Q=1$). In both plots we see that NPV is indeed inversely correlated with the sum of the squared velocities, and that the maximum NPV corresponds to ${\textbf f}^*$. These results support the choice of the objective function used in the first ST optimization step (Eq.~\ref{well_rate_ratios_opt_eqn}). It thus appears that, even though many physical aspects of the problem are not included in the first ST step, the method can nonetheless provide useful results for well-rate ratios, at least for this example.

It is also noteworthy that the ${\textbf f}^*$ solution provides the maximum NPV with $d=0.1$ (Fig.~\ref{fig:npv_vs_uTu}b). This suggests that minimizing the squared velocities is a reasonable treatment even with some amount of discounting. We believe there are cases, however, for which this direct correspondence will no longer hold. For example, in a problem with large $d$ (e.g., $d \gtrsim 0.2$), we would expect early production to be very highly valued, and this effect might not be captured in the first ST step. Other economic aspects, such as very high injected water cost, might also represent important effects not captured in the first ST step. It is possible, however, that some of these effects could be (at least partially) treated by modifying the first-step optimization problem, or by adjusting the ${\textbf W}_d$ matrices in Eq.~\ref{well_rate_ratios_opt_eqn}.
\color{black}

Next, in order to assess the degree of improvement achieved by the two-step ST (and by MADS and adjoint-gradient-based optimization), we evaluate NPV for three base cases for each of the well scenarios shown in Fig.~\ref{fig:100x100_3example_cases}. The base cases entail equally distributed injection and production rates; i.e., $f_k^p=1/n_p$ for $k=1, \ldots, n_p$ and $f_k^i=1/n_i$ for $k=1, \ldots, n_i$, for field rates that cover a range of PVI values. We ensure that the range includes the $Q^*$ value provided by ST.

\begin{table}[!htb]
    \centering
     \caption{NPVs for base cases ($M=3$) for 2D reservoir model}
    \begin{tabular}{*4c}
	   \hline
      {}         & \multicolumn{3}{c}{}	\\
      Well       & \multicolumn{3}{c}{NPV (BN \$)}		\\
      scenario   & \multicolumn{3}{c} {  }   \\
       {}        & PVI = 0.25 & 0.75 & 1.25   \\
		\hline
     $1^1$    &      0.54 & 0.48 & 0.15   	\\
     $2^2$    &      0.81 & 0.89 & 0.62    	\\
     $3^3$    &      0.49 & 1.22 & 1.41    	\\
       \hline     
       \multicolumn{4}{l} {$\textnormal{}^1{Q}^{*}_{1, M=3}=0.49~\textnormal{PVI}$} \\
       \multicolumn{4}{l} {$\textnormal{}^2{Q}^{*}_{2, M=3}=0.41~\textnormal{PVI}$} \\
       \multicolumn{4}{l} {$\textnormal{}^3{Q}^{*}_{3, M=3}=0.65~\textnormal{PVI}$} \\
	\end{tabular}
    \label{tab:st2D_base_NPVs}
 \end{table}
 
The base-case NPVs are presented in Table~\ref{tab:st2D_base_NPVs}. The $Q^*$ value for each case is also given.  Although the use of equally distributed rates may seem somewhat simplistic, we did confirm that, for the best base cases, the water-oil-ratio (WOR) in the production wells does not reach the point where water production cost exceeds oil revenue. This means that the use of a standard reactive control procedure (in which production wells are shut in if they become uneconomic) in these cases would not have any impact on NPV. Thus, the best base cases are reasonable in this regard. 
 
The base case results can then be compared to those in Table~\ref{tab:st2D_NPVs}, where we show optimization results for each configuration for three different mobility ratios. This table includes ST results along with those using MADS and adjoint gradients. Results for the latter two procedures are median NPVs from five separate runs. The minimum and maximum NPVs for these methods are shown in Table~\ref{tab:st2D__max-min_NPVs}. We reiterate that ST provides (essentially) global optima in both steps, so it does not need to be run multiple times to avoid poor local optima.
 
Base case and ST results are also compared, in terms of NPV versus time plots for each well configuration (for $M=3$), in Fig.~\ref{fig:100x100_3example_cases}d,e,f. These results, along with those in Tables~\ref{tab:st2D_base_NPVs} and \ref{tab:st2D_NPVs}, demonstrate the clear improvement provided by our two-step ST. This improvement is due both to the optimization of well-rate ratios and to the optimization of $Q$, with the relative impact of these two steps varying from case to case. Of essential interest is the close correspondence between ST optimization results and those from MADS and adjoint gradients (Table~\ref{tab:st2D_NPVs}). The largest relative discrepancy is for well scenario~2, $M=5$, where we see ST underperform MADS by 5.6\%. The close correspondence between results from our two-step ST and those from formal optimization methods is very encouraging, and suggests that the surrogate treatments may indeed capture important aspects of the true full-physics optimization problem.
 
The results in Table~\ref{tab:st2D__max-min_NPVs} demonstrate that MADS achieves very similar results in all five runs for each case. This is also generally observed for the adjoint-gradient procedure, though the well scenario~3, $M=5$ case is an exception.  For the minimum NPV (1.11~BN~\$) run, we applied a tighter stopping criterion and performed additional SNOPT iterations, but the NPV did not improve. Thus, this result may correspond to convergence to a relatively poor local optimum.
\color{black}

\begin{table}[!htb]
    \centering
     \caption{Optimal NPVs using different methods (2D reservoir model)}
    \begin{tabular}{*5c}
	   \hline
      {}   &        \multicolumn{4}{c}{}	\\
      Well     & Mobility   &    \multicolumn{3}{c}{NPV (BN \$)}		\\
      scenario & ratio      &  {}   &   {}	    &   Adjoint-  	\\
       {}   & {}   & ST    &   $\textnormal{MADS}^1$	&   $\textnormal{gradient}^1$  	\\
		\hline
     {}   &      1	  &   1.23    &   1.24	&   1.24   	\\
     1    &      3	  &   0.62    &   0.62	&   0.62  	\\
     {}   &      5	  &   0.47    &   0.47	&   0.47   	\\
       \hline       
     {}   &      1	  &   2.03    &   2.02	&   2.07   	\\
     2    &      3	  &   1.03    &   1.09	&   1.02   	\\
     {}   &      5	  &   0.84    &   0.89	&   0.89   	\\
       \hline     
     {}   &      1	 &   2.55    &   2.64	&   2.55    	\\
     3    &      3	 &   1.47    &   1.47	&   1.47  	\\
     {}   &      5	 & 1.28      &   1.29	&   1.29   	\\
       \hline
     \multicolumn{5}{l} {$\textnormal{}^1$Median values}
	\end{tabular}
    \label{tab:st2D_NPVs}
\end{table}

\begin{table}[!htb]
    \centering
     \caption{Minimum and maximum NPVs for MADS and adjoint-gradients (2D reservoir model)}
    \begin{tabular}{*7c}
	   \hline
      \multicolumn{7}{c}{}	\\
      Well     & Mobility&  \multicolumn{4}{c}{NPV (BN \$)}		\\
      scenario & ratio   &  {}   &  {}   &   {}   &   \multicolumn{2}{c} {Adjoint-}  	\\
       {}      &   {}    & ST    & \multicolumn{2}{c}{MADS} & \multicolumn{2}{c}{gradient}  \\
       {}      &   {}    & {}    & min  & max   & min  & max  \\
		\hline
     {}   &      1	 & 1.23 &   1.23    &   1.24	&   1.24  & 1.24	\\
     1    &      3	 & 0.62 &   0.61    &   0.62	&   0.62  & 0.62	\\
     {}   &      5	 & 0.47 &   0.46    &   0.47	&   0.46  & 0.47 	\\
       \hline       
     {}   &      1	 & 2.03 &   2.02    &   2.02	&   2.07  & 2.07 	\\
     2    &      3	 & 1.03 &   1.09    &   1.09	&   1.02  & 1.03 	\\
     {}   &      5	 & 0.84 &   0.89    &   0.89	&   0.87  & 0.89	\\
       \hline     
     {}   &      1	& 2.55 &   2.64  &   2.64	&   2.55  & 2.55 \\
     3    &      3	& 1.47 &   1.47  &   1.47	&   1.43  & 1.47 \\
     {}   &      5	& 1.28 &   1.28  &   1.29	&   1.11  & 1.29 \\
       \hline     
	\end{tabular}
    \label{tab:st2D__max-min_NPVs}
\end{table}

Elapsed (wall-clock) times for the optimizations are shown in Table~\ref{tab:st2D_elapsed_time}. These values are median timings for the five MADS and adjoint-gradient runs (recall ST need only be run once). Individual full-physics simulation runs for these cases require $\sim$40--50 seconds. We see that ST, as noted earlier, entails an elapsed time corresponding to 1--2 full-physics runs. MADS requires five or more iterations for these cases  (we verified that additional iteration with tighter termination criteria did not improve the MADS results), so the speedup using ST is fairly significant (speedup factor of $6\times$ or more). Speedup relative to adjoint-gradient-based optimization is larger because this procedure does not readily parallelize, in contrast to MADS.

\begin{table}[!htb]
    \centering
    \caption{Elapsed time (fully parallelized for ST and MADS) for different methods for 2D reservoir model}
    \begin{tabular}{*5c}
	   \hline
      {}   &        \multicolumn{4}{c}{}	\\
       Well    &     Mobility     &    \multicolumn{3}{c}{Elapsed time (s)}		\\
      scenario &   ratio       &  {}      &   {}	    &   Adjoint-  	\\
       {}   &    {}       &  ST    &   $\textnormal{MADS}^1$	&   $\textnormal{gradient}^1$   	\\
		\hline
     {}   &      1	 &   52    &  389   &  1652    	\\
     1    &      3	 &   52    &  385   &  1952    	\\
     {}   &      5	 &   51    &  389   &  3665     	\\
       \hline       
     {}   &      1	 &   60    &   361   &  1138   	\\
     2    &      3   &   63    &   458   &  1165   	\\
     {}   &      5	 &   62    &   456   &  2538   	\\
       \hline     
     {}   &      1	&   67     &   503   &   876  	    \\
     3    &      3	&   68     &   503   &  2750    	\\
     {}   &      5	&   66     &   502   &  3690   	\\
     \hline 
     \multicolumn{5}{l} {$\textnormal{}^1$Median values} \\     
	\end{tabular}
    \label{tab:st2D_elapsed_time}
\end{table}

The results in this section demonstrate that our two-step surrogate treatment provides reasonable accuracy relative to formal optimization methods, and that it can lead to significant speedup. As such, it may be useful for both standalone well control optimization applications, and as an inner-loop in field development optimization. The latter point, highlighted by the ST results for many different configurations (Fig.~\ref{fig:npv_st_vs_mads}), will be exploited in a future study. We next consider the application of ST for a 3D example.

\subsection{3D reservoir model}
\label{subsec:4.2}

\begin{figure*}
  \includegraphics[width=.2\textwidth]{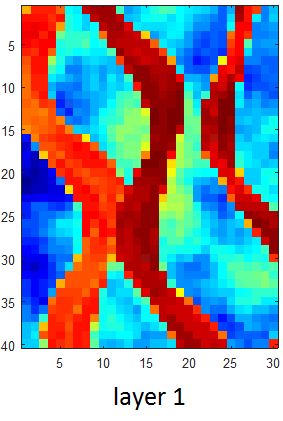}\hfill
  \includegraphics[width=.2\textwidth]{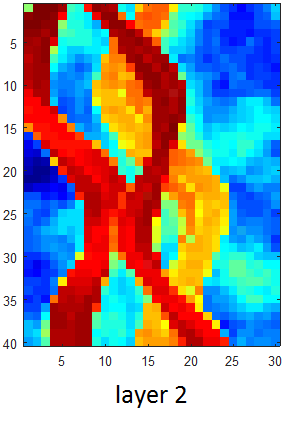}\hfill
  \includegraphics[width=.2\textwidth]{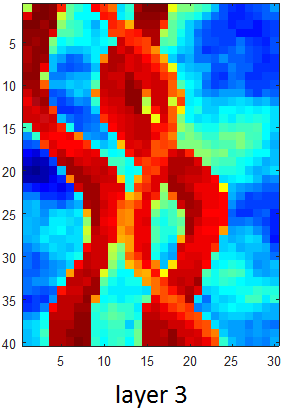}\hfill
  \includegraphics[width=.2\textwidth]{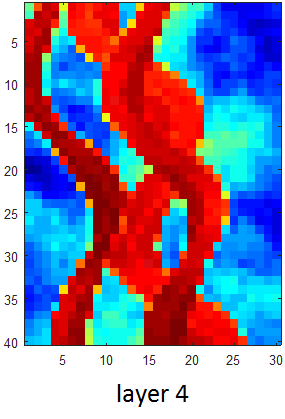}\hfill
  \includegraphics[width=.2\textwidth]{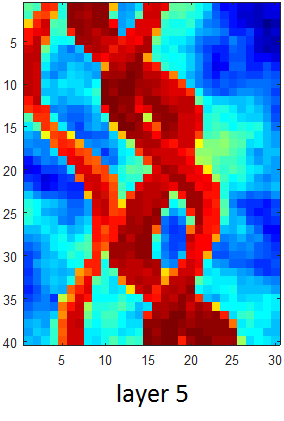}
  \includegraphics[width=.2\textwidth]{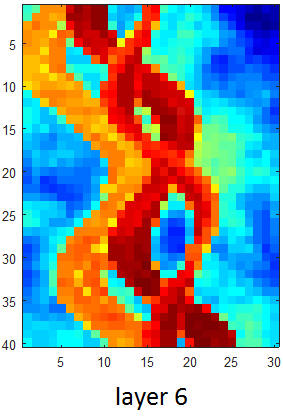}\hfill
  \includegraphics[width=.2\textwidth]{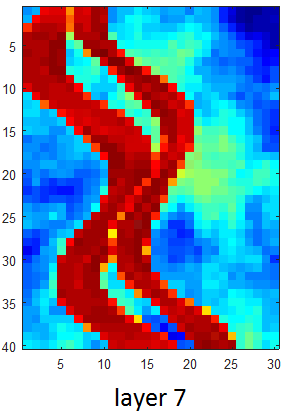}\hfill
  \includegraphics[width=.2\textwidth]{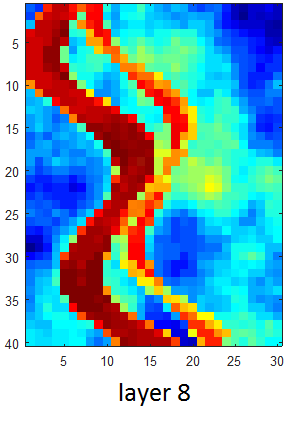}\hfill
  \includegraphics[width=.2\textwidth]{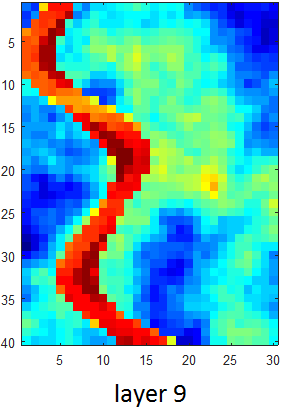}\hfill
  \includegraphics[width=.2\textwidth]{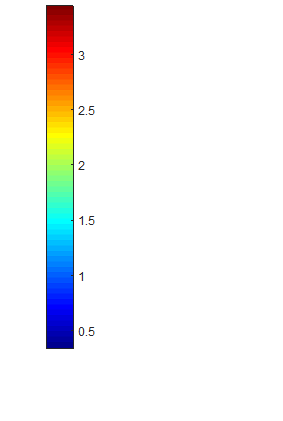}
  \caption{Layer-by-layer permeability field (in the $x$-direction) for the $40 \times 30 \times 9$ model representing a prograding fluvial channel system (model from \cite{castro}).}
  \label{fig:40x30x9_permeability_x}
\end{figure*}

\begin{figure*}[t!]
    \centering
    \begin{subfigure}[b]{.5\textwidth}
        \centering
        \includegraphics[width=\textwidth]{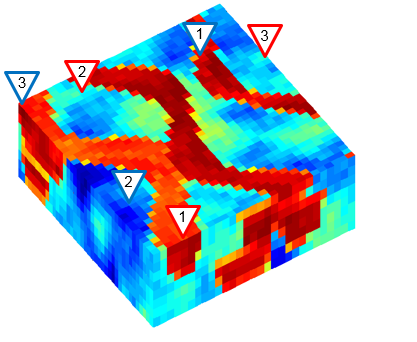}
        \caption{Well scenario 1}
    \end{subfigure}%
    ~ 
    \begin{subfigure}[b]{.5\textwidth}
        \centering
        \includegraphics[width=\textwidth]{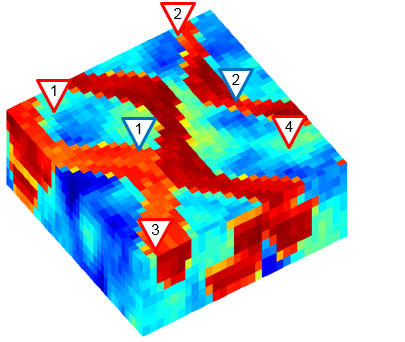}
        \caption{Well scenario 2}
    \end{subfigure}   
    \caption{Well configurations considered for the 3D case.}
    \label{fig:40x30x9_2well_scenario}
\end{figure*}

The 3D model used in this assessment is shown in Fig.~\ref{fig:40x30x9_permeability_x}. This $40 \times 30 \times 9$ model (10,800 total grid blocks) represents a prograding fluvial channel system, and corresponds to a portion of a synthetic geological model developed by Castro~\cite{castro}. We assess ST performance for two different scenarios involving (1) three producers and three injectors, and (2) four producers and two injectors. All wells are completed over the full reservoir thickness. The two scenarios are shown in Fig.~\ref{fig:40x30x9_2well_scenario}. Both single and multiple control periods are considered for this case.

Relevant simulation and optimization parameters appear in Table~\ref{tab:40x30x9_simul_opt_pars}. We again consider $M=1$, 3 and 5, and use the same relative permeability curves as in the 2D cases (the curves in Fig.~\ref{fig:krel_curve} for $M\neq 1$, and straight-line relative permeabilities for $M=1$).

\begin{table}[!htb]
    \centering
    \caption{Simulation and optimization parameters for $40 \times 30 \times 9$ channelized model}
    \begin{tabular}{ll}
		\hline
		Grid dimensions                     & 40  $\times$ 30  $\times$ 9 \\
        Grid cell dimensions                & 100 ft $\times$ 100 ft $\times$ 15 ft \\
		Initial pressure $p_i$, at datum    & 6000 psi at 6000 ft\\
		$\mu_o$						        & 1.0, 3.0 or 5.0 cp \\
        $\mu_w$ (constant)			        & 1.0 cp \\
		$\rho_o$ and $\rho_w$               & 53.10 and 64.79 lbm/ft$^3$ \\
		$B_o$ and $B_w$ at $p_i$            & 1.075 RB/STB \\
		Oil compressibility at $p_i$		& 10$^{-5}$ psi$^{-1}$ \\
		Water compressibility at $p_i$		& 10$^{-5}$ psi$^{-1}$ \\
        Rock compressibility at $p_i$		& 10$^{-9}$ psi$^{-1}$  \\
        Field porosity						& 0.25 \\
        Simulation period                   & 14,600 days \\
        $p_o$, $c_{pw}$ and $c_{iw}$        & \$60, \$5 and \$5/STB \\
        $c_w$								& \$20 MM \\
        Discount rate                       & 0\% or 5\% \\
        \hline
	\end{tabular}
    \label{tab:40x30x9_simul_opt_pars}
\end{table}


\subsubsection{Single control period}
\label{subsec:4.2.1}

Optimization results for the two scenarios shown in Fig.~\ref{fig:40x30x9_2well_scenario} will now be presented. We first consider the case where well rates are fixed for the entire simulation time frame (as was the treatment in all of the 2D runs). As in the 2D case, several base cases, with equal distribution of rates among all injection and all production wells, are assessed for $M=3$.  Again, for the best base cases, none of the production wells reach the point where water production cost exceeds oil revenue. As noted earlier, this indicates that the use of standard reactive controls would not impact the best base-case NPVs. \color{black}

The base-case NPVs are shown in Table~\ref{tab:3D_base_NPVs}, and the optimized NPVs using ST and MADS are presented in Table~\ref{tab:opt3D_NPVs}. The median, minimum and maximum NPVs using MADS are provided.
Adjoint-gradient-based optimizations were not performed for this case, though we would expect these results to closely coincide with the MADS results, as in the 2D example. 

\begin{figure*}[t!]
    \centering
    \begin{subfigure}[b]{0.32\textwidth}
        \centering
         \includegraphics[width=\textwidth]{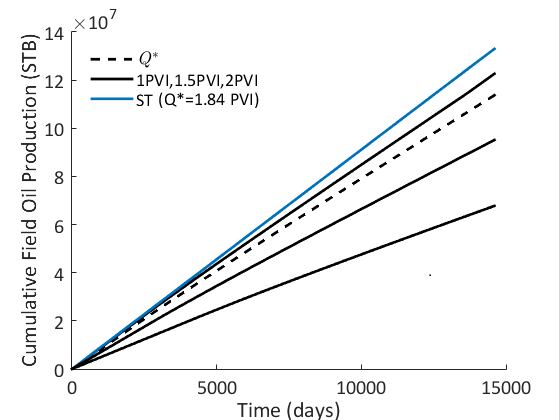}
        \caption{Cumulative field oil production}
    \end{subfigure}%
    ~ 
    \begin{subfigure}[b]{0.32\textwidth}
        \centering
        \includegraphics[width=\textwidth]{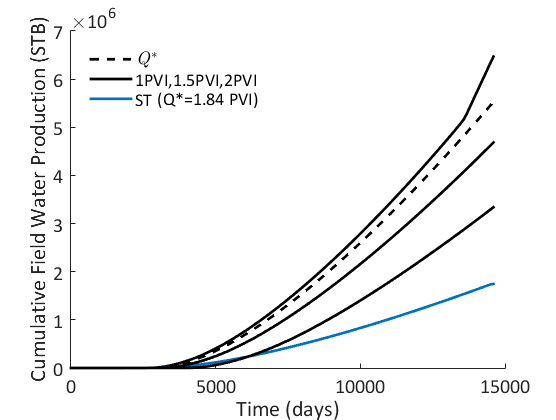}
        \caption{Cumulative field water production}
    \end{subfigure}   
    ~ 
    \begin{subfigure}[b]{0.32\textwidth}
        \centering
         \includegraphics[width=\textwidth]{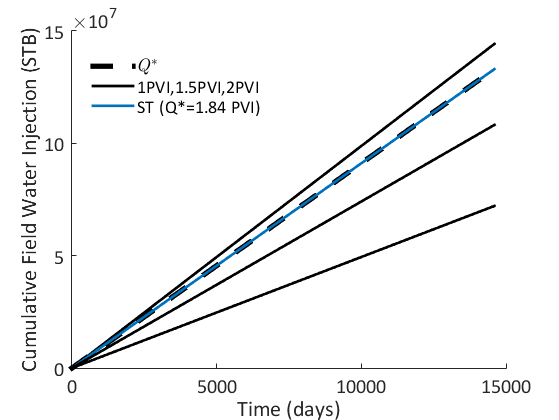}
        \caption{Cumulative field water injection}
    \end{subfigure}   
    \caption{Cumulative field oil production (a), water production (b), and water injection (c) for base cases ($Q^*$, 1~PVI, 1.5~PVI and 2~PVI) and ST (single control period, well scenario~2, $M=3$, 3D example).}
    \label{fig:cum_prod_inj_base_st_well_scenario2}
\end{figure*}

As in the 2D cases, we again see that the NPVs provided by ST clearly exceed those for the three base cases for both scenarios. Base-case and ST results for cumulative field oil and water production and cumulative water injection, for well scenario~2 with $M=3$, are shown in Fig.~\ref{fig:cum_prod_inj_base_st_well_scenario2}. Results are also presented for total injection (and production) of $Q^*$, as determined by ST, but with equal injection and production between wells.  The ST results display the most oil production, along with the least water production, of the cases considered. These results highlight the importance of optimizing both ${\textbf{f}}$ and $Q$.

\begin{table}[!htb]
    \centering
     \caption{NPVs for base cases ($M=3$) for 3D reservoir model with a single control period}
    \begin{tabular}{*4c}
	   \hline
      {}       &        \multicolumn{3}{c}{}	\\
      Well     &         \multicolumn{3}{c}{NPV (BN \$)}		\\
      scenario & \multicolumn{3}{c} { }  \\
       {}      &    PVI = 1.0 & 1.5 & 2.0   \\
		\hline
     $1^1$         &    3.71 & 4.72 & 5.22  	\\
         
     $2^2$         &    3.69 & 5.14 & 6.61   	\\
     
       \hline     
       \multicolumn{4}{l} {$\textnormal{}^1{Q}^{*}_{1, M=3}=1.44~\textnormal{PVI}$} \\
       \multicolumn{4}{l} {$\textnormal{}^2{Q}^{*}_{2, M=3}=1.84~\textnormal{PVI}$} \\
	\end{tabular}
    \label{tab:3D_base_NPVs}
 \end{table}
 
\begin{table}[!htb]
    \centering
    \caption{Optimal NPVs for ST and MADS for 3D reservoir model with a single control period}
    \begin{tabular}{*6c}
	   \hline
      {}   &        \multicolumn{3}{c}{}	\\
      Well     &     Mobility   &    \multicolumn{4}{c}{NPV (BN \$)}		\\
      scenario &     ratio      &  {}   &   \multicolumn{3}{c}{MADS}	     	\\
        {}  &   {}    & ST &   med  & min & max	  	\\
		\hline
     {}   &      1	 &    17.64 &   17.44	&   17.35    &  17.53  	\\
     1   &      3	 &    5.45  &   5.52	&    5.48    &   5.60  	\\
     {}   &      5	 &    3.85  &   4.04	&    4.04    &   4.04  	\\
       \hline       
     {}   &      1	 &  22.70 &   22.88	 &   22.88    &   22.88	 	\\
     2    &      3	 &  7.32  &   7.56	 &    7.56    &   7.57	 	\\
     {}   &      5	 &  4.65  &   4.80 	&     4.80    &   4.80	\\
       \hline     
       \end{tabular}
    \label{tab:opt3D_NPVs}
\end{table}

Comparison of ST and MADS results in Table~\ref{tab:opt3D_NPVs} demonstrates that ST again provides results quite close to those achieved using a formal optimization method. The largest relative difference between ST and (median) MADS results is for well scenario~1, $M=5$, where ST underperforms MADS by 4.7\%. For well scenario~1, $M=1$, ST actually achieves a larger NPV than any of the MADS runs. We also see that there is very little variation between the MADS results in most cases. 
 For well scenario~1 with $M=1$ and 3, the MADS results do vary by 1--2\% between the minimum and maximum NPV runs. For the minimum-NPV runs, we tightened the termination criterion and performed additional MADS iterations, but this did not provide improvement in the NPVs. This observation supports the idea that the different local optima may correspond to slightly different NPVs. \color{black}

Table~\ref{tab:rmsv3D_elapsed_time} presents elapsed time, under full parallelization, for ST and MADS optimizations. Individual simulation runs require $\sim$105--115 seconds for these cases. We thus see that ST elapsed times continue to be less than the time required for two full-physics runs. MADS requires six iterations for these cases, and speedup using ST of about $5\times$ relative to MADS is observed. It is again apparent that ST is able to achieve reasonable results along with substantial speedup. 

\begin{table}[!htb]
    \centering
    \caption{Elapsed time (fully parallelized) for ST and MADS for a single control period for 3D reservoir model}
    \begin{tabular}{*4c}
	   \hline
      {}   &  {} &  \multicolumn{2}{c}{}	\\
       Well    &     Mobility     &    \multicolumn{2}{c}{Elapsed time (s)}		\\
      scenario &   ratio       &  {}      &   {}	     	\\
       {}   &    {}  &  ST    &   $\textnormal{MADS}^1$		\\
		\hline
     {}   &      1	 &  151    &   741   	\\
     1    &      3	 &  146    &   736	   	\\
     {}   &      5	 &  147    &   738   	\\
       \hline       
     {}   &      1	 &  132    &   741   	\\
     2    &      3   &  136    &   735   	\\
     {}   &      5	 &  137    &   734   \\
       \hline      
       \multicolumn{4}{l} {$\textnormal{}^1$Median values}
    \end{tabular}
    \label{tab:rmsv3D_elapsed_time}
\end{table}


\subsubsection{Multiple control periods}
\label{subsec:4.2.2}

Our examples up to this point have involved a single control period and zero discount rate (except for the $M=1$ results in Fig.~\ref{fig:npv_vs_uTu}b, where we considered $d=0.1$). We now relax these assumptions and consider multiple control period strategies, with and without discounting. Both cases involve well scenario~1 (three producers and three injectors). In the first case we consider $M=3$ and $d=0$, and in the second case $M=5$ and $d=0.05$. The well rates are now determined four times during the simulation time frame, which increases $n_{\textrm{opt}}$ to 20 (16 well-rate ratios plus 4 field rates). Two approaches can be used to assess ST performance relative to MADS for this case. First, we can simply apply the ST well rates over the entire time frame. Second, the single-period ST solution can be used as an initial guess for MADS (we designate this case ST \& MADS).

\begin{figure}
  \includegraphics[width=.5\textwidth]{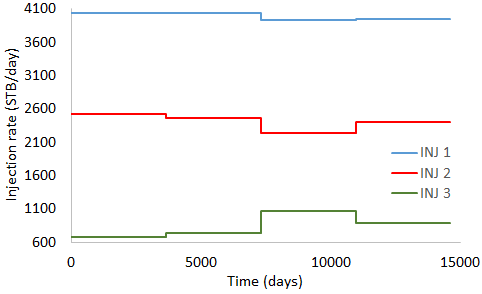}
  \caption{Optimized injection rates obtained by ST \& MADS (3D model, multiple control periods, well scenario~1, $M=3$, $d=0$).}
  \label{fig:opt_inj_rates}
\end{figure}

\begin{figure}
  \includegraphics[width=.5\textwidth]{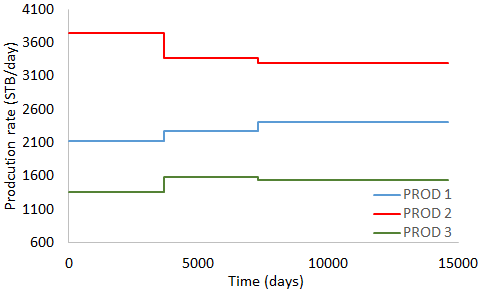}
  \caption{Optimized production rates obtained by ST \& MADS (3D model, multiple control periods, well scenario~1, $M=3$, $d=0$).}
  \label{fig:opt_prod_rates}
\end{figure}

The NPV results and corresponding elapsed times obtained with ST (standard single control period solution), ST \& MADS, and MADS (we give the median, minimum and maximum NPVs for the five standalone MADS runs) are presented in Tables~\ref{tab:st_3D_NPVs_multiple} and \ref{tab:st_3D_time_multiple}. We see that both methods provide generally similar NPV results. For the first case, there does not appear to be substantial benefit from using multiple control periods (compare results in Tables~\ref{tab:opt3D_NPVs} and \ref{tab:st_3D_NPVs_multiple}). We do, however, observe that the NPV resulting from the ST \& MADS approach provides an objective function value that is 3\% higher than the median solution from standalone MADS for the first case.  Even though this difference is relatively small, this result suggests that the ST solution is indeed a reasonable initial guess for the subsequent MADS run.  The corresponding optimized injection and production rates for the four control intervals, obtained by ST \& MADS for the $M=3$, $d=0$ case, are presented in Figs.~\ref{fig:opt_inj_rates} and~\ref{fig:opt_prod_rates}. We do see some shifts from control period to control period, though these shifts are relatively slight. 

For the second case ($M=5$, $d=0.05$), the ST result is 2.8\% below the median MADS result, but it is above the minimum MADS result. The ST \& MADS approach provides the same objective function value as the median MADS case. These results, along with those in Fig.~\ref{fig:npv_vs_uTu}b, suggest that our two-step ST is indeed applicable for cases with nonzero discount rate. As noted earlier, the first ST step may become less reliable in cases with large $d$, where very early production is highly beneficial. The effect of large $d$ would, however, be captured in the second ST step.
\color{black}

\begin{table}[!htb]
    \centering
    \caption{Optimal NPVs for 3D reservoir model with multiple control periods (well scenario~1).}
    \begin{tabular}{*7c}
	   \hline
          {} & {} &  \multicolumn{5}{c}{NPV (BN \$)}	\\
        Mobility & d  & ST  &  ST \& &   \multicolumn{3}{c}{MADS}	\\
         ratio   & (\%) & {}  &  MADS &   med & min & max 	\\
          \multicolumn{7}{c}{}	\\
		\hline
            3 & 0 & 5.45   &  5.63 & 5.46  & 5.36  & 5.55	\\     
            5 & 5 & 1.71   &  1.76 & 1.76  & 1.46  & 1.76 	\\
       \hline       
	\end{tabular}
    \label{tab:st_3D_NPVs_multiple}
\end{table}

\begin{table}[!htb]
    \centering
    \caption{Elapsed time (fully parallelized) for 3D reservoir model with multiple control periods (well scenario~1).}
    \begin{tabular}{*6c}
	   \hline
          {} &  \multicolumn{5}{c}{Elapsed time (s)}	\\
        Mobility &   ST   &  ST \& &   \multicolumn{3}{c}{MADS}	\\
         ratio &   {}   &  MADS &   med & min & max 	\\
          \multicolumn{6}{c}{}	\\
		\hline
            3 &  136   &  971  & 1355 & 1228 & 1377 	\\     
            5 &  137   &  1211 & 1515 & 1493 & 1655 	\\
       \hline       
	\end{tabular}
    \label{tab:st_3D_time_multiple}
\end{table}

\begin{figure}
  \includegraphics[width=.4\textwidth]{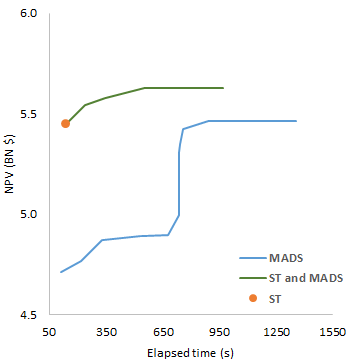}
  \caption{Evolution of NPV for ST, ST \& MADS, and standalone MADS (3D model, multiple control periods, well scenario~1, $M=3$, $d=0$).}
  \label{fig:elapsed_time_vs_npv_M3_mult}
\end{figure}

From the timings in Table~\ref{tab:st_3D_time_multiple}, we see that ST \& MADS is slightly faster than MADS (ST \& MADS requires 70--80\% of the median time for MADS). The improvement in NPV with iteration for the first case ($M=3$, $d=0$) is shown in Fig.~\ref{fig:elapsed_time_vs_npv_M3_mult}, which displays the evolution of NPV for ST, ST \& MADS, and standalone MADS (median case). In this case, ST \& MADS leads to speedup relative to standalone MADS because it requires 9~iterations, while standalone MADS performs 13~ iterations (both methods use the same termination criteria). 

A potential advantage of ST \& MADS is that standalone MADS might be run multiple times to find the best solution (even though the differences in NPVs between MADS runs are typically rather small), while ST \& MADS needs to be run only once. This is possible because the initial guess provided by standard ST is expected to lead to a promising (local) optimum, as is indeed observed for this example. This behavior, if achieved consistently over a range of more challenging test cases, would represent a real advantage for ST \& MADS, as multiple optimization runs would not be required.

\section{Concluding remarks}
\label{sec:5}

In this work, we introduced a new two-step surrogate treatment (ST) for well control optimization of oil fields under waterflood. The method provides injection and production rates for all wells for a single control period. Our ST proceeds by first determining  appropriate well-rate ratios, which are expected to lead to more uniform reservoir sweep. In the second step, overall field rate is optimized such that the objective function of the actual problem (NPV in the cases considered here) is maximized or minimized. In the first step we use simplified physics, specifically the assumption of unit-mobility-ratio displacement and exact voidage replacement. This step is accelerated through use of velocity response matrices, which can be constructed efficiently in a preprocessing step. Full-physics simulations are used in the second ST step, but this minimization can be performed quickly since it involves only a single optimization variable. Near-global optima are obtained in both steps, though this does not provide the global optimum to the actual problem due to the approximations introduced in the ST. The overall two-step procedure requires an elapsed time corresponding to the runtime for 1--2 full-order simulations, assuming full parallelization in the second ST step.

The ST was applied for both 2D and 3D optimization problems. For the 2D case, 40 well scenarios (corresponding to a wide range of well locations and types), generated by PSO during a field development optimization run, were used to compare ST and MADS. Relatively close correspondence between the two sets of optimization results, for $M=3$, was observed. More detailed comparisons were presented for three of the well scenarios for systems with $M=1$, 3 and 5. Multiple MADS and adjoint-gradient runs were performed to assure avoidance of poor local optima with these methods  (though this is generally not a major issue with oil-water well control optimization problems of the type considered here).  Again, close correspondence between results from our two-step ST and the formal optimization procedures was observed. Similar observations applied for the 3D case, where two well scenarios were considered, again for $M=1$, 3 and 5. In terms of timings, under full parallelization ST provided speedups of $5 \times$ or more relative to MADS for the cases considered. For a case involving multiple control periods, where the ST result provided the initial guess for MADS, a speedup of less than $2 \times$ was observed. However, for such cases standalone MADS would likely be run multiple times, while the ST \& MADS approach need be run only once.

There are a number of areas in which future work extending our two-step ST could be performed. More general well control optimization setups, involving, e.g., additional cases with nonzero discount rates or the consideration of optimal economic project life~\cite{epl}, could be considered. Optimization under geological uncertainty is another important aspect that should be addressed. In this case realization selection, using procedures described in \cite{representative}, could be applied. The two-step ST may be particularly useful within the context of field development optimization, where well number, locations, types and controls must be determined. Such MINLP problems are extremely demanding computationally, so surrogate treatments can have a large impact. The applicability of the two-step ST for this problem is suggested by the results shown in Fig.~\ref{fig:npv_st_vs_mads}, and our findings in a followup study indeed confirm this. The two-step ST could also be extended to handle more complex production scenarios, such as primary production followed by waterflood, or water-alternating-gas processes.  Further testing with realistic 3D models should also be performed. Finally, comparisons with other heuristic procedures for determining well allocation/rates (such as those in \citep{torrado,ushmaev}) should be performed. Optimization results using the various approaches could then be assessed, and these findings could be used to devise a hybrid implementation based on the most promising components from the different methodologies. \color{black}

\section*{Acknowledgements}

We thank Nikhil Padhye for providing the Unified Optimization Framework (which facilitated testing and comparing the different optimization algorithms used in this study) and for helping us with its use. We are also grateful to Marco Antonio Cardoso for useful discussions and suggestions. We thank Petrobras and the industrial affiliates of the Stanford Smart Fields Consortium for partial funding of this work.

\section*{References}

\bibliography{mybibfile}

\end{document}